\documentclass[12pt]{amsart}
\advance\hoffset by -0.5 truecm
\usepackage{amssymb}
\newtheorem{Theorem}{Theorem}[section]
\newtheorem{Lemma}[Theorem]{Lemma}
\newtheorem{Corollary}[Theorem]{Corollary}

\newtheorem{Proposition}[Theorem]{Proposition}

\newtheorem{Definition}[Theorem]{Definition}
\newtheorem{Example}[Theorem]{Example}
\newtheorem{Remark}[Theorem]{Remark}

\def \dim{{\mbox {dim}}\,}

\def\V{\mbox{Var}}
\newcommand{\comp}{\mbox{\tiny{o}}}

\def\Z{{\mathbb Z}}
\def\R\re
\def\V{\bf V}

\def \re{{\mathbb R}}
\def \Q{{\mathbb Q}}
\def \cp{{\mathbb CP}}
\def \T{{\mathbb T}}
\def \C{{\mathbb C}}

\def \0{\lambda_{0}}

\def\h{{\rm h}_{\rm top}(g)}
\def\T{{\mathcal T}}
\def\en{{\rm h}_{\rm top}}

\begin{document}
\title[Minimal entropy and collapsing]{Minimal entropy and
collapsing with curvature bounded from below}

\author[G. P. Paternain]{Gabriel P. Paternain}\thanks{G.P. Paternain is on leave from Centro de Matem\'atica, Facultad de Ciencias, Igu\'a 4225, 11400 Montevideo, Uruguay.}
\address{CIMAT  \\
          A.P. 402, 36000 \\
          Guanajuato. Gto. \\
          M\'exico.}
\email {paternain@cimat.mx}

\author[J. Petean]{Jimmy Petean}
 \address{CIMAT  \\
          A.P. 402, 36000 \\
          Guanajuato. Gto. \\
          M\'exico.}
\email{jimmy@cimat.mx}

\thanks{J. Petean is supported by grant 28491-E of CONACYT}


\date{2000}


\begin{abstract} We show that if a closed manifold $M$ admits an
$\mathcal{F}$-structure (not necessarily polarized, possibly of
rank zero) then its minimal entropy vanishes. In particular, this
is the case if $M$ admits a non-trivial $S^1$-action. As a
corollary we obtain that the simplicial volume of a manifold
admitting an $\mathcal{F}$-structure is zero.

We also show that if $M$ admits an ${\mathcal F}$-structure then
it collapses with curvature bounded from below. This in turn
implies that $M$ collapses with bounded scalar curvature or,
equivalently, its Yamabe invariant is non-negative.

We show that ${\mathcal F}$-structures of rank zero appear rather
frequently: every compact complex elliptic surface admits one as
well as any simply connected closed 5-manifold.

We use these results to study the minimal entropy problem. We show
the following two theorems: suppose that $M$ is a closed manifold
obtained by taking connected sums of copies of $S^{4}$, $\cp^{2}$,
$\overline{\cp}^{2}$, $S^{2}\times S^{2}$ and the $K3$ surface.
Then $M$ has zero minimal entropy. Moreover, $M$ admits a smooth
Riemannian metric with zero topological entropy if and only if $M$
is diffeomorphic to $S^{4}$, $\cp^{2}$, $S^{2}\times S^{2}$,
$\cp^{2}\#\overline{\cp}^{2}$ or $\cp^{2}\#\cp^{2}$. Finally,
suppose that $M$ is a closed simply connected 5-manifold. Then $M$
has zero minimal entropy. Moreover, $M$ admits a smooth Riemannian
metric with zero topological entropy if and only if $M$ is
diffeomorphic to $S^{5}$, $S^{3}\times S^{2}$, the nontrivial
$S^{3}$-bundle over $S^{2}$ or the Wu-manifold $SU(3)/SO(3)$.

\end{abstract}

\maketitle

\section{Introduction}

Let $M^{n}$ be a closed orientable connected smooth manifold.
Given a Riemannian metric $g$, let $\phi_{t}$ be the geodesic flow
of $g$.

Perhaps the simplest dynamical invariant that one can associate to
$\phi_{t}$ to roughly measure its orbit structure complexity is
the {\it topological entropy,} which we denote by $\h$. Positive
entropy means in general, that the geodesic flow presents
somewhere in the phase space (the unit sphere bundle of the
manifold) a complicated dynamical behaviour. There are various
equivalent ways of defining entropy (see Subsection \ref{entcurv})
and among them there is a formula, known as Ma\~n\'e's formula,
that gives a nice Riemannian description of $\h$. Given points $p$
and $q$ in $M$ and $T>0$, define $n_{T}(p,q)$ to be the number of
geodesic arcs joining $p$ and $q$ with length $\leq T$. R.
Ma\~n\'e \cite{Man} showed that $$\h= \lim_{T\rightarrow
\infty}\frac{1}{T}\log \int_{M\times M}n_{T}(p,q)\;dp\,dq. $$

One of the main goals in this paper will be the study of the
variational theory of the functional $g\mapsto \h$. In general
this functional is only upper semicontinuous in the $C^{\infty}$
topology (\cite{N,Y}) and it has a simple behaviour under scaling
of the metric: if $c$ is any positive constant, then ${\rm
h}_{{\rm top}}(cg)=\frac{\h}{\sqrt{c}}$. Hence if we want to
extract interesting extremal metrics from this functional a
normalization is required. The Riemannian invariant that we will
use for this normalization is the volume ${\rm Vol}(M,g)$.

Set the {\it minimal entropy} of $M$ to be
\[{\rm h}(M):=\inf\{{\rm h}_{\rm top}(g)\: |\: g \mbox{ is a smooth
metric on $M$ with } {\rm Vol}(M,g) =1\}. \] A smooth metric
$g_{0}$ with ${\rm Vol}(M, g_0) =1$ is {\it entropy minimizing} if
\[{\rm h}_{\rm top}(g_{0})={\rm h}(M). \]

The {\it minimal entropy problem} for $M$ is whether or not there
exists an entropy minimizing metric on $M$.  Say that the {\it
minimal entropy problem can be solved} for $M$ if there exists an
entropy minimizing metric on $M$. Smooth manifolds are hence
naturally divided into two classes: those for which the minimal
entropy problem can be solved and those for which it cannot.
Passing by, we note that we do not obtain a meaningful invariant
if we replace the infimum by the supremum. Indeed, Manning proved
in \cite{Ma2} that
\[\sup\{{\rm h}_{\rm top}(g)\: |\: g \mbox{ is a smooth
metric on $M$ with } {\rm Vol}(M,g) =1\}=\infty. \]

There are a number of classes of manifolds for which the minimal
entropy problem can be solved.  For instance, the minimal entropy
problem can always be solved for a closed orientable surface $M$.
For the 2-sphere and the 2-torus, this follows from the fact that
both a metric with constant positive curvature and a flat metric
have zero topological entropy.  For surfaces of higher genus, A.
Katok \cite{Ka} proved that each metric of constant negative
curvature minimizes topological entropy, and conversely that any
metric that minimizes topological entropy has constant negative
curvature.

This result of Katok has been generalized to higher dimensions by
G. Besson, G. Courtois and S. Gallot \cite{BCG}, as follows.
Suppose that $M^{n}$ ($n\geq 3$) admits a locally symmetric metric
$g_{0}$ of negative curvature, normalized so that ${\rm
Vol}(M,g_0) =1$. Then $g_{0}$ is the unique entropy minimizing
metric up to isometry. Unlike the case of a surface, a locally
symmetric metric of negative curvature on a closed $n$-manifold
($n\ge 3$) is unique up to isometry, by the rigidity theorem of
Mostow.

A positive solution to the minimal entropy problem appears to
single out manifolds that have either a high degree of symmetry or
a low topological complexity. What this means in our context will
become apparent below.  A similar phenomena is observed for closed
3-manifolds \cite{AP}.

There is a close relationship between minimal entropy, {\it
minimal volume} and {\it simplicial volume}. As we shall explain
in Subsection \ref{chain} there is a positive constant $c(n)$ such
that
\begin{equation}
c(n)\,\|M\|\leq [{\rm h}(M)]^{n}\leq (n-1)^{n}\mbox{\rm
MinVol}(M). \label{schain}
\end{equation}
Recall that the minimal volume MinVol($M$) is the infimum of ${\rm
Vol}(M,g)$ where $g$ runs over all metrics whose sectional
curvature is bounded in absolute value by 1. Also recall that the
simplicial volume of a closed orientable manifold $M$, $\|M\|$, is
defined as the infimum of $\sum_{i}|r_{i}|$ where the $r_{i}$ are
the coefficients of a real cycle that represents the fundamental
class of $M$. This number is a homotopy invariant of $M$. The
minimal volume does depend on the smooth structure of $M$ (see
\cite{Bessi}) but we do not know if the same holds true for the
minimal entropy.

Computing these invariant is in general a very difficult task. J.
Cheeger and M. Gromov introduced in \cite{Cheeger, Gromov} the
concept of ${\mathcal F}$-structure (see Section \ref{F-str} for
the precise definition). An {\it ${\mathcal F}$-structure} on a
manifold $M$ is a natural generalization of an effective torus
action on $M$. The structure partitions $M$ into disjoint orbits
which are flat manifolds amenable to collapse. When the dimension
of the orbits are, in a certain precise sense, locally constant,
the structure is said to be {\it polarized}. The simplest
${\mathcal F}$-structures are the $\T$-structures, which consist
of open coverings of the manifold and torus actions on each of the
elements of the covering which commute on overlaps. The simplest
example of a polarized $\T$-structure is given by a locally free
circle action.

Cheeger and Gromov proved in \cite{Cheeger, Gromov} that if $M$
admits a polarized ${\mathcal F}$-structure then the minimal
volume of $M$ vanishes. The vanishing of the minimal volume
implies in turn that all the characteristic numbers of the
manifold are zero. Cheeger and Gromov also proved that if the
${\mathcal F}$-structure has {\it positive rank}, i.e., all its
orbits have positive dimension, then the Euler characteristic of
$M$ must be zero. There exist plenty of examples of closed
manifolds which admit $\mathcal{F}$-structures but whose Euler
characteristic is non-zero. Therefore they do not admit
$\mathcal{F}$-structures of positive rank. For instance the Euler
characteristic of any simply connected closed 4-manifold is
strictly positive, but for any $m,n,k,l$, the manifold $n\cp^2 \#
k{\overline \cp}^2\# mK3 \# l(S^2 \times S^2)$ admits a
$\T$-structure of rank zero. This will follow form the results in
Section 5. Hence, general ${\mathcal F}$-structures are abundant
in comparison with polarized ones.

 In Section \ref{zerome}
we show:

\medskip
\noindent {\bf Theorem A.} {\it If the closed  manifold $M$ admits
an $\mathcal{F}$-structure then the minimal entropy of $M$ is 0. }

\medskip

The theorem and (\ref{schain}) yield the following corollary,
which generalizes the result of K. Yano \cite{Yano} that closed
manifolds which admit non-trivial $S^1$-actions have simplicial
volume 0 (there is also a proof of the latter result in
\cite{Gromov}).

\medskip
\noindent{\bf Corollary.} {\it Let $M$ be a closed manifold. If
$M$ admits an $\mathcal{F}$-structure then the simplicial volume
of $M$ is 0.}
\medskip

Hence the existence of an ${\mathcal F}$-structure, possibly of
rank zero, also imposes constraints on the topology of the
manifold.

\medskip

The method employed in the proof of Theorem A is general enough
that allows us to apply it to the study of other types of
collapsing. We will say that $M$ {\it collapses with curvature
bounded from below} if there exists a sequence of metrics $g_{j}$
for which the sectional curvature is uniformly bounded from below,
but their volumes approach zero as $j$ goes to infinity. Similarly
we will say that $M$ collapses with Ricci (respectively, scalar)
curvature bounded from below if there exists a sequence of metrics
$g_{j}$ for which the Ricci (respectively, scalar) curvature is
uniformly bounded from below, but their volumes approach zero as
$j$ goes to infinity. In Section \ref{yamabe} we prove:

\medskip
\noindent {\bf Theorem B.} {\it If the closed manifold $M$ admits
an $\mathcal{F}$-structure then $M$ collapses with curvature
bounded from below. }

\medskip

Clearly if $M$ collapses with curvature bounded form below then it
also collapses with Ricci and scalar curvatures bounded from
below. As we explain in Section \ref{yamabe} if $M$ has dimension
$\geq 3$, then it collapses with scalar curvature bounded form
below if and only if it collapses with {\it bounded scalar
curvature}. This is in turn equivalent to having non-negative {\it
Yamabe invariant}.

It is interesting to remark that for instance the manifold
$T^{4}\#{\overline\cp}^{2}$ admits an ${\mathcal F}$-structure but
it does {\it not} collapse with bounded Ricci curvature (see
Section \ref{yamabe}). Therefore our Theorem B can be regarded as
an optimal extension of the results of Cheeger and Gromov in the
sense that there is no stronger collapsing phenomena for general
${\mathcal F}$-structures other than the one claimed in the
theorem.

C. LeBrun proved in \cite{LeBrun2} that 
the Yamabe invariant of any compact complex
surface of general type is strictly negative. It follows
from Theorem B that these surfaces do not admit ${\mathcal
F}$-structures. Among these surfaces of general type there are
simply connected ones which are homeomorphic (but not
diffeomorphic) to connected sums of $\cp^{2}$'s and
$\overline{\cp}^{2}$'s. Hence in dimension 4 there are simply
connected closed manifolds which do not admit ${\mathcal
F}$-structures and they are homeomorphic to manifolds that do
admit them. We do not know if this a phenomena specific of
dimension 4. In dimension $\geq 5$ the second author showed in
\cite{Petean} that any simply connected manifold has non-negative
Yamabe invariant. This opens the possibility that {\it any} closed
simply connected manifold of dimension $\geq 5$ admits an
${\mathcal F}$-structure. Morever it is possible for this
structure to be polarized in odd dimensions. In fact we show in
Section \ref{me45}:

\medskip
\noindent {\bf Theorem C.} {\it Every simply connected closed
smooth 5-manifold $M$ admits a ${\mathcal T}$-structure. Moreover,
suppose that either:

\begin{enumerate}
\item $M$ is spin;
\item $M$ is the non-trivial $S^{3}$-bundle over $S^{2}$ or the Wu-manifold $SU(3)/SO(3)$;
\item $M$ is a connected sum of manifolds of types 1 or 2.
\end{enumerate}
Then $M$ admits a polarized ${\mathcal T}$-structure.}
\medskip

We do not know if {\it every} closed simply connected non-spin
5-manifold admits a polarized ${\mathcal T}$-structure, even
though it appears to be the case.

These results can be used to give fairly complete solutions to the
minimal entropy problem for simply connected manifolds of
dimensions 4 and 5.

\medskip
\noindent {\bf Theorem D.} {\it Let $M$ be a closed manifold
obtained by taking connected sums of copies of $S^{4}$, $\cp^{2}$,
$\overline{\cp}^{2}$, $S^{2}\times S^{2}$ and the $K3$ surface.
Then ${\rm h}(M)=0$ and the minimal entropy problem can be solved
for $M$ if and only if $M$ is diffeomorphic to $S^{4}$, $\cp^{2}$,
$S^{2}\times S^{2}$, $\cp^{2}\#\overline{\cp}^{2}$ or
$\cp^{2}\#\cp^{2}$.}

\medskip

A manifold $M$ like in Theorem D realizes many intersection forms
of simply connected 4-manifolds. In fact the 11/8-conjecture
(see \cite{Do1, Do2}) states that any smooth simply connected 4-manifold is
homeomorphic to a manifold as in Theorem D. Hence, if one assumes
the 11/8-conjecture, Theorem D is saying that any smooth simply
connected 4-manifold is homeomorphic to one whose minimal entropy
is zero and for which we know the answer to the minimal entropy
problem.

The proof of Theorem D is partially based on the fact that the
$K3$ surface admits a $\T$-structure. In fact we show that any
elliptic compact complex surface 
admits a $\T$-structure. We also show that $\T$-structures behave
relatively well with respect to the usual operations of connected
sums and surgeries on manifolds.

For simply connected 5-manifolds, we have a complete answer to the
minimal entropy problem :

\medskip

\noindent {\bf Theorem E.} {\it Let $M$ be a closed simply
connected 5-manifold. Then ${\rm h}(M)=0$ and the minimal entropy
problem can be solved for $M$ if and only if $M$ is diffeomorphic
to $S^{5}$, $S^{3}\times S^{2}$, the nontrivial $S^{3}$-bundle
over $S^{2}$ or the Wu-manifold $SU(3)/SO(3)$. }

\medskip

The common feature of the nine manifolds listed in Theorems D and
E is that they are {\it elliptic}. This means that their loop
space homology grows polynomially for every coefficient field (cf.
Section 3, \cite{FH1,FH2,GH1} and references therein). In fact, as
we will see in Section 3, these are the only elliptic manifolds in
dimensions 4 and 5. Hence Theorems D and E characterize this very
much studied class of manifolds as that for which the minimal
entropy problem can be solved or, equivalently, as that for which
there exists a smooth metric $g$ with $\h=0$. It is tempting to
speculate that perhaps the same phenomena occurs in any dimension.

We would like to close this introduction by illustrating some of
the ideas with specific examples. A 5-dimensional Brieskorn
variety of type $(a_{1},a_{2},a_{3},a_{4})$ is given by the
intersection of the 7-sphere in $\C^{4}$ with the zero set of:
\[f(z_{1},z_{2},z_{3},z_{4})=z_{1}^{a_{1}}+z_{2}^{a_{2}}+z_{3}^{a_{3}}
+z_{4}^{a_{4}}.\] This gives a large class of simply connected
5-manifolds. In fact, they are all spin and their second homology
group can be computed using the algorithm described in \cite{O,R}.
The Brieskorn varieties admit very simple polarized ${\mathcal
T}$-structures: they have a canonically defined action of $S^{1}$
which is locally free (but not free in general). If we let
$$q_{i}=\mbox{\rm lcm}(a_{1},a_{2},a_{3},a_{4})/a_{i}$$ then the
action is given by:
\[e^{i\theta}(z_{1},z_{2},z_{3},z_{4})=
(e^{q_{1}i\theta}z_{1},e^{q_{2}i\theta}z_{2},e^{q_{3}i\theta}z_{3},e^{q_{4}i\theta}z_{4}).\]
For example, the Brieskorn variety $M_{2}$ defined by $(2,3,3,3)$
coincides with the spin manifold whose second homology group is
$\Z_{2}\oplus\Z_{2}$. It has the property that its loop space
homology grows exponentially with $\Z_{2}$ coefficients and hence
for any $C^{\infty}$ Riemannian metric $g$, $\h>0$ (see Theorem
\ref{poshtop} in Section \ref{me45}). Since ${\rm
MinVol}(M_{2})={\rm h}(M_{2})=0$ it follows that the minimal
entropy problem for $M_{2}$ cannot be solved. It is interesting to
note that $M_{2}$ has the rational cohomology ring of the 5-sphere
and hence its loop space homology with rational coefficients is
actually bounded, i.e., $M_{2}$ is rationally elliptic.

 \vspace{.5cm}

{\it Acknowledgement:} We thank S. Halperin for explaining us how
to compute the growth of the loop space homology of a 5-manifold.

\section{Preliminaries on simplicial volume, minimal volume and topological entropy}
\label{introduction}

The purpose of this Section is to present some of the basic
material and definitions that we will need later on.

\subsection{Simplicial volume}

Let $M$ be a closed manifold.  Denote by $C_{*}$ the real chain
complex of $M$: a chain $c\in C_{*}$ is a finite linear combination
$\sum_{i}r_{i}\sigma_{i}$ of singular simplices $\sigma_{i}$ in $M$
with real coefficients $r_{i}$.  Define the {\it simplicial
$l^{1}$-norm} in $C_{*}$ by setting $|c|=\sum_{i}|r_{i}|$. This norm
gives rise to a pseudo-norm on the homology $H_{*}(M,\re)$ by setting
\[| [\alpha] |=\inf\{ |z|\: :\: z\in C_* \mbox{ and } [z] =[\alpha]
\}. \]
When $M$ is orientable, define the {\it simplicial volume} of $M$,
denoted $\|M\|$, to be the simplicial norm of the fundamental class.
The simplicial volume is also called {\it Gromov's invariant}, since
it was first introduced by Gromov in \cite{Gromov}.

\subsection{Minimal volume and collapsing}

The {\it minimal volume}\index{minimal volume} MinVol($M$) of a
manifold $M$ is defined to be the infimum of $\mbox{\rm Vol}(M,g)$
over all metrics $g$ in ${\mathcal R}(M)$ such that the sectional curvature
$K_{g}$ of $g$ satisfies $|K_{g}|\leq 1$. This differential invariant
was introduced by M. Gromov in \cite{Gromov}.

As we mentioned in the introduction we have \cite{Cheeger}:

\begin{Proposition} If $M$ admits a polarized ${\mathcal F}$-structure, then $\mbox{\rm MinVol}(M)=0$.
\label{minvol0}
\end{Proposition}

\subsection{Topological entropy and curvature}
\label{entcurv} We recall in this subsection the definition of the
topological entropy of the geodesic flow of a Riemannian metric
$g$ on a closed manifold $M$. The geodesic flow of $g$ is a flow
$\phi_{t}$ that acts on $SM$, the unit sphere bundle of $M$, which
is a closed hypersurface of the tangent bundle of $M$. In general
the topological entropy is defined for an arbitrary continuous
flow (or map) on a compact metric space.

Let $(X,d)$ be a compact metric space and let $\phi_{t}:X
\rightarrow X$ be a continuous flow. For each $T>0$ we define a
new distance function
\[d_{T}(x,y):= \max_{0\leq t\leq T}\,d(\phi_{t}(x),\phi_{t}(y)).\]
Since $X$ is compact, we can consider the minimal number of balls
of radius $\varepsilon>0$ in the metric $d_{T}$ that are necessary
to cover $X$. Let us denote this number by $N(\varepsilon,T)$. We
define
\[{\rm h}(\phi,\varepsilon):= \limsup_{T\rightarrow \infty}\frac{1}{T}\log N(\varepsilon,T).\]
Observe now that the function
 $\varepsilon\mapsto {\rm h}(\phi,\varepsilon)$ is monotone decreasing
and therefore the following limit
 exists:
\[{\rm h}_{{\rm top}}(\phi):= \lim_{\varepsilon\rightarrow 0}h(\phi,\varepsilon).\]
The number $\en(\phi)$ thus defined is called the {\it topological
entropy} of the flow $\phi_{t}$. Intuitively, this number measures
of orbit complexity of the flow. The positivity of $\en(\phi)$
indicates complexity or ``chaos'' of some kind in the dynamics of
$\phi_{t}$. The topological entropy $\en(\phi)$ may also be
defined as $\en(\phi_{1})$ using the entropy of the time one-map
or it may be defined in either of the following ways. All the
definitions give the same number $\en(\phi)$ which is independent
of the choice of metric \cite{HK,W}.

A set $Y\subset X$ is called a $(T,\varepsilon)$-separated set if
given different points $y,y'\in Y$ we have
$d_{T}(y,y')\geq\varepsilon$. Let $S(T,\varepsilon)$ denote the
maximal cardinality of a $(T,\varepsilon)$-separated set.
 Then
\[\en(\phi)=\lim_{\varepsilon\to 0}\limsup_{T\rightarrow\infty}\frac{1}{T}\log S(T,\varepsilon).\]

A set $Z\subset X$ is called a $(T,\varepsilon)$-spanning set if
for all $x\in X$ there exists $z\in Z$ such that
$d_{T}(x,z)\leq\varepsilon$.
 Let $M(T,\varepsilon)$ denote the minimal cardinality of a $(T,\varepsilon)$-spanning set.
Then
\[\en(\phi)=\lim_{\varepsilon\to 0}\limsup_{T\rightarrow\infty}\frac{1}{T}\log M(T,\varepsilon).\]

Given a compact subset $K\subset X$ (not necessarily invariant) we
can define the topological entropy of the flow with respect to the
set $K$, $\en(\phi,K)$, simply by considering
 separated (spanning) sets of $K$.

The following proposition gives an idea of the dynamical
significance of the topological entropy (for proofs see
\cite{HK,W}).

\begin{Proposition}The topological entropy verifies the following properties:
\begin{enumerate}
\item For any two closed subsets $Y_{1}$, $Y_{2}$ in $X$,
\[\en(\phi,Y_{1}\cup Y_{2})=\max_{i=1,2}\,\en(\phi,Y_{i});\]

\item If $Y_{1}\subset Y_{2}$ then $\en(\phi,Y_{1})
\leq \en(\phi,Y_{2})$;

\item Let $\phi_{t}^{i}:X_{i} \rightarrow X_{i}$ for $i=1,2$ be two flows and let $\pi:X_{1}\rightarrow X_{2}$ be a continuous map commuting with $\phi_{t}^{i}$ i.e. $\phi_{t}^{2}\comp\pi=\pi\comp \phi_{t}^{1}$.
 If $\pi$ is onto, then $\en(\phi^{1})\geq \en(\phi^{2})$ and if $\pi$ is finite-to-one, then $\en(\phi^{1})\leq \en(\phi^{2})$.

\item Let $\phi_{t}^{i}:X_{i} \rightarrow X_{i}$ for $i=1,2$ be two flows
and let $\psi_{t}:=\phi_{t}^{1}\times \phi_{t}^{2}$ be the product
flow on $X_{1}\times X_{2}$. Then
$\en(\psi)=\en(\phi^{1})+\en(\phi^{2})$.

\item Given $c\in\re$, let $c\phi_{t}$ be the flow given by
$c\phi_{t}:=\phi_{ct}$. Then $\en(c\phi)=|c|\en(\phi)$.

\end{enumerate}

\label{propiedades-3}
\end{Proposition}

Next we shall state a useful result of R. Bowen that we will need
later.

\begin{Proposition}[Corollary 18 in \cite{Bo}] Let $(X,d)$ and $(Y,e)$ be compact metric spaces and $\phi_{t}:X\rightarrow X$ a flow. Suppose $\pi:X\rightarrow Y$ is a continuous map such that $\pi\comp \phi_{t}=\pi$.
 Then
\[\en(\phi)=\sup_{y\in Y}\,\en(\phi,\pi^{-1}(y)).\]
\label{integral1-3}
\end{Proposition}

Given a Riemannian metric $g$, let $d$ be any distance function
compatible with the topology of $SM$. Since the geodesic flow is a
smooth flow on $SM$ we can attach to it its topological entropy
that we denote by $\h$ to stress its dependence on the Riemannian
metric $g$. There is a formula, known as Ma\~n\'e's formula, that
gives a nice alternative way of thinking about $\h$. Given $p$ and
$q$ in $M$ and $T>0$, define $n_{T}(p,q)$ as the number of
geodesic arcs joining $p$ and $q$ with length $\leq T$. R.
Ma\~n\'e showed in \cite{Man} that
  $$\h=
  \lim_{T\rightarrow \infty}\frac{1}{T}\log
  \int_{M\times M}n_{T}(p,q)\;dp\,dq.
   $$

Using property 5 in Proposition \ref{propiedades-3} it is easy to
check how entropy behaves under scaling: if $c$ is any positive
constant, then ${\rm h}_{{\rm top}}(cg)=\frac{\h}{\sqrt{c}}$.

We now describe a basic relationship between entropy and
curvature.

Let $(M^{n},g)$ be a closed Riemannian manifold and let
$K_{max}$ be a positive upper bound for the sectional curvature. It was proved in
\cite{Paternain} that
\[{\rm h}_{\rm top}(g)\leq \frac{n-1}{2}\sqrt{K_{max}}-\frac{\min_{v\in SM}\,r(v)}{2\sqrt{K_{max}}},\]
where $SM$ is the unit sphere bundle of $M$ and $r(v)$ is the Ricci curvature
in the direction of $v\in SM$.

Let $k$ be a positive number such that $|K(P)|\leq k$ for all 2-planes $P$. Then, clearly
$r\geq -(n-1)k\,g$ and hence the previous inequality gives
\begin{equation}
{\rm h}_{\rm top}(g)\leq \frac{n-1}{2}\sqrt{k}+\frac{n-1}{2}\sqrt{k}=(n-1)\sqrt{k}.
\label{manni}
\end{equation}
The latter inequality was first proved by A. Manning in \cite{Ma2}.

\subsection{An important chain of inequalities}
\label{chain}
 Let $(M,g)$ be a closed Riemannian manifold and let
$\widetilde{M}$ be its universal covering endowed with the induced
metric.  Given $x\in\widetilde{M}$, let $V(x,r)$ be the volume of
the ball with center $x$ and radius $r$.  Set
\[\lambda(g):=\lim_{r\to +\infty}\frac{1}{r}\log\,V(x,r).\]
Manning \cite{Ma} showed that the limit exists and it is independent
of $x$.

Set
\[\lambda(M):=\inf\{\lambda(g)\: |\: g \mbox{ is a smooth metric on
$M$ with } {\rm Vol}(M,g) =1\}. \]

It is well known \cite{Mil} that $\lambda(g)$ is positive if and
only if $\pi_{1}(M)$ has exponential growth. Manning's inequality
\cite{Ma} asserts that for any metric $g$,
\begin{equation}
\lambda(g)\leq {\rm h}_{\rm top}(g).
\label{manineq}
\end{equation}
In particular, it follows that if $\pi_{1}(M)$ has exponential
growth then ${\rm h}_{\rm top}(g)$ is positive for any metric $g$.
This fact was first observed by E.I. Dinaburg in \cite{D}.  Gromov
showed in \cite{Gromov} that if ${\rm Vol}(M,g) = 1$, then there
is a positive constant $c(n)$ such that
\begin{equation}
c(n)\|M\|\leq [\lambda(g)]^n. \label{simvol-5}
\end{equation}
Finally it was observed in \cite{P} that using (\ref{manni}) it is
easy to show that
\begin{equation}
[{\rm h}(M)]^{n}\leq (n-1)^{n}\mbox{\rm MinVol}(M).
\label{pineq}
\end{equation}
Hence if we combine (\ref{manineq}), (\ref{simvol-5}) and
(\ref{pineq}), we obtain the following chain of inequalities:
\begin{equation}
c(n)\|M\|\leq  [\lambda(M)]^{n} \leq [{\rm h}(M)]^{n}\leq
(n-1)^{n}\mbox{\rm MinVol}(M). \label{chain-5}
\end{equation}
The only known manifolds with ${\rm h}(M)>0$ are manifolds with $\|M\|\neq
0$. For these manifolds $\pi_{1}(M)$ has exponential growth.

\subsection{Entropy of products and submersions}

\begin{Lemma}
\begin{enumerate}
\item Let $(M_{1},g_{1})$ and $(M_{2},g_{2})$ be two compact Riemannian
manifolds. Endow $M_{1}\times M_{2}$ with the product metric
$g_{1}\times g_{2}$. Then
$${\rm h}_{\rm top}(g_{1}\times g_{2})=\sqrt{[{\rm h}_{\rm top}(g_{1})]^{2}+[{\rm h}_{\rm top}(g_{2})]^{2}}.$$
\item Let $(M,g_{M})\mapsto (N,g_{N})$ be a Riemannian submersion
where $M$ and $N$ are compact manifolds.
Then ${\rm h}_{\rm top}(g_{M})\geq {\rm h}_{\rm top}(g_{N})$.
\end{enumerate}
\label{facil}
\end{Lemma}

\begin{proof}Let us prove the first item.
Let $f:S(M_{1}\times M_{2})\to S^{1}$ be the function given by
\[f(x_{1},v_{1},x_{2},v_{2})=(|v_{1}|_{x_{1}},|v_{2}|_{x_{2}}).\]
Since the geodesics in $M_{1} \times M_{2}$ are products of
geodesics in $M_{1}$ and $M_{2}$, the function $f$ is constant
along the orbits of the geodesic flow of $M_{1}\times M_{2}$. It
follows from Proposition \ref{integral1-3} that
\[{\rm h}_{\rm top}(g_{1}\times g_{2})=\sup_{c\in S^{1}}{\rm h}_{top}(f^{-1}(c)).\]
If we write $c=(l,m)$, it is easy to check using Proposition
\ref{propiedades-3} that $${\rm h}_{\rm top}(f^{-1}(c))=l\,{\rm
h}_{\rm top}(g_{1})+m\,{\rm h}_{\rm top}(g_{2})$$ from which we
obtain right away the first equality in the lemma.

To prove the second item, let $H\subset SM$ be the set of all
horizontal unit vectors. Clearly the geodesic flow of $(M,g_{M})$
leaves $H$ invariant. Let $\tau:H\to SN$ be the restriction to $H$
of the differential of the submersion map. Since horizontal
geodesics project to geodesics, $\tau$ is a surjective map that
intertwines the geodesic flow of $(M,g_{M})$ restricted to $H$
with the geodesic flow of $(N,g_{N})$. It follows from Proposition
\ref{propiedades-3} that ${\rm h}_{\rm top}(g_{M})\geq {\rm
h}_{\rm top}(g_{N})$.

\end{proof}

\section{Elliptic manifolds in dimensions 4 and 5}

Let $M$ be a closed simply connected manifold and let $\Omega M$
be the space of based loops. Let $k_{p}$ be the prime field of characteristic
$p$, $p$ prime or zero.
Following Y. F\'elix, S. Halperin
and J.C. Thomas we say that $M$ is {\it elliptic} if for each $p$, the homology of the loop space:
\[\sum_{i=0}^{n}\mbox{\rm dim}\,H_{i}(\Omega M,k_{p}),\]
grows polynomially with $n$ (cf. \cite{FH1,FH2,GH1} and references therein).

Elliptic manifolds are rare. However a number of geometrically
interesting spaces are elliptic:

\begin{enumerate}

\item homogeneous spaces;

\item manifolds $M$ admitting a fibration
$F\to M\to B$ with $F$ and $B$ elliptic;

\item manifolds $M$ for which the algebra $H^{*}(M,k_{p})$
is generated by two elements for all $p$;

\item manifolds $M$ admitting a smooth action by a compact Lie group
with a simply connected codimension one orbit;

\item connected sums $M\# N$ with the algebras $H^{*}(M,\Z)$
and $H^{*}(N,\Z)$ each generated by a single class.

\end{enumerate}

The manifold $M$ is said to be {\it rationally elliptic}
if the total rational homotopy $\pi_{*}(M)\otimes\Q$ is finite
 dimensional, i.e.
there exists a positive integer $i_{0}$ such that for all $i\geq
i_{0}$, $\pi_{i}(M)\otimes\Q=0$. This property is known to be
equivalent to the polynomial growth of $\sum_{i=0}^{n}\mbox{\rm
dim}\,H_{i}(\Omega M,\Q)$. Obviously an elliptic manifold is
rationally elliptic. We will see that for smooth 4-manifolds
ellipticity and rational ellipticity are equivalent. This is no
longer the case for 5-manifolds as we will see below.

\begin{Lemma}Suppose that $M$ is 4-dimensional and let $b_{2}$ be the
second Betti number of $M$. If $M$ is rationally elliptic then
$b_{2}\leq 2$.  \label{lema1}
\end{Lemma}

\begin{proof} It is shown in \cite[Corollary 1.3]{FrH} (cf. also \cite{Fel}) that
if $M^{n}$ is rationally elliptic then,
\begin{equation}
\sum_{k\geq 1}2k\,\dim\,(\pi_{2k}(M)\otimes\Q)\leq n.
\label{frh}
\end{equation}
Since $M$ is simply connected the Hurewicz isomorphism theorem
implies that
\[b_{2}=\dim\,H_{2}(M,\Q)=\dim\,(\pi_{2}(M)\otimes\Q).\]
Since $n=4$, using (\ref{frh}) we obtain $2\,b_{2}\leq 4$.

\end{proof}

\begin{Lemma}Let $M$ be a closed smooth simply connected
4-manifold. The following are equivalent:
\begin{enumerate}
\item $M$ is elliptic;
\item $M$ is rationally elliptic;
\item $M$ is homeomorphic to $S^{4}$, $\cp^{2}$,
$S^{2}\times S^{2}$, $\cp^{2}\#\overline{\cp}^{2}$ or
$\cp^{2}\#\cp^{2}$.
\end{enumerate}
Moreover, if $M$ is not elliptic then $\sum_{i=0}^{n}\mbox{\rm
dim}\,H_{i}(\Omega M,\Q)$ grows exponentially.
\label{lema2}
\end{Lemma}

\begin{proof} Obviously 1 implies 2.
Let us prove that 2 implies 3. Suppose that $M$ is rationally elliptic. By Lemma \ref{lema1},
$b_{2}\leq 2$. Since $M$ is smooth, the Kirby-Siebenmann obstruction
vanishes. Therefore by M. Freedman's theory \cite{Fr}, the homeomorphism
type of $M$ is completely determined by the intersection form of $M$.
It follows that if $b_{2}=0$, $M$ is homeomorphic to $S^{4}$ and if
$b_{2}=1$, $M$ is homeomorphic to $\cp^{2}$.
When $b_{2}=2$, the possible intersection forms are
\[\left(\begin{array}{cr}0&1\\1&0\end{array}\right),\;\;\;\;
\left(\begin{array}{cr}1&0\\0&-1\end{array}\right)\;\;\;\mbox{\rm and}\;\;\;
\left(\begin{array}{cr}1&0\\0&1\end{array}\right).\]
These forms correspond to $S^{2}\times S^{2}$, $\cp^{2}\#\overline{\cp}^{2}$ and $\cp^{2}\#\cp^{2}$ respectively.

On the other hand $S^{4}$, $\cp^{2}$ and $S^{2}\times S^{2}$ are
homogeneous spaces and hence they are elliptic (see property 1 above).
By property 5 above,
$\cp^{2}\#\overline{\cp}^{2}$ and $\cp^{2}\#\cp^{2}$ are
elliptic.

Finally, it is well known that the homology of the loop space with
rational coefficients can either grow polynomially or
exponentially.
\end{proof}

\begin{Remark}{\rm For an arbitrary simply connected manifold
$M$ it is known that if $\sum_{i=0}^{n}\mbox{\rm
dim}\,H_{i}(\Omega M,k_{p})$ does not grow polynomially then it
must grow at least like $\lambda^{\sqrt{n}}$ for some $\lambda>1$
\cite{FH2}. There is a conjecture that says that the growth should
in fact be exponential, but this is only known for rational
coefficients (as we mentioned at the end of the proof of the last
lemma) and for primes $p$ strictly bigger that the dimension of
$M$.}
\end{Remark}

\begin{Theorem}[Following a suggestion of S. Halperin] Let $M$ be a closed
$(2s-1)$-connected manifold of dimension $4s+1$ with $s\geq 1$.
Then $M$ is elliptic if and only if $H_{2s}(M,\Z)$ is $0$, $\Z$ or
$\Z_{2}$. Moreover, if $M$ is not elliptic the homology of the
loop space of $M$ grows exponentially for some field of
coefficients $k_{p}$. \label{loop5m}
\end{Theorem}

\begin{proof}It follows from a theorem of S. Eilenberg and J.C. Moore \cite[Theorem 12.1]{EM}
that the homology of the loop space can be computed as
\[H_{*}(\Omega M,k_{p})\cong\mbox{\rm Tor}^{C^{*}(M)}(k_{p},k_{p}),\]
where $C^{*}(M)$ is the differential graded algebra given by the
normalized singular cochains with coefficients in $k_{p}$. (In
fact, Eilenberg and Moore mention in his paper that this special
case of their theorem has to be attributed to J.F. Adams
\cite{Ad}.)

It can be seen that for a manifold $M$ satisfying the hypotheses
of the theorem there exists a {\it quism} between $C^{*}(M)$ and
$(H^{*}(M,k_{p}),0)$. This means a morphism of differential graded
algebras with the property that induces isomorphisms in homology.
Since a quism preserves $\mbox{\rm Tor}$ it follows that
\[H_{*}(\Omega M,k_{p})\cong\mbox{\rm Tor}^{H^{*}(M,k_{p})}(k_{p},k_{p}).\]
We now make use of the following lemma whose proof will be given after completing the proof
of the theorem.

\begin{Lemma}The sum of the dimensions of $\mbox{\rm Tor}^{H^{*}(M,k_{p})}(k_{p},k_{p})$ grows exponentially
unless $\mbox{\rm dim}\,H_{2s}(M,k_{p})\leq 1$. Conversely if
$\mbox{\rm dim}\,H_{2s}(M,k_{p})\leq 1$ then the sum of the
dimensions of $\mbox{\rm Tor}^{H^{*}(M,k_{p})}(k_{p},k_{p})$ grow
polynomially.
\end{Lemma}

A result of C.T.C. Wall \cite{Wa} (see also the corollary before
Lemma F in \cite{Barden}) using the linking form ensures that the
torsion part of $H_{2s}(M,\Z)$ always has the form $B+B$ or
$B+B+\Z_{2}$ for some finite abelian group $B$. Hence if $M$ is
elliptic, the lemma implies that $B$ must be zero and when the
$\Z_{2}$ factor appears the rank of $H_{2s}(M,\Z)$ should be zero.
\qedsymbol

\medskip

\noindent{\it Proof of the lemma.} Let us set for brevity
$k:=k_{p}$. Observe that $R:=H^{*}(M,k)$ is a (graded) commutative
local ring with residue field $k$ that satisfies Poincar\'e
duality. We note that it suffices to prove the lemma ignoring the
grading of $R$ because $\mbox{\rm Tor}_{p,q}^{R}(k,k)=0$ for
$q>p(4s+1)$ (the first integer indicates the resolution degree and
the second the internal grading).

Let $a:=\mbox{\rm dim}\,H_{2s}(M,k)=\mbox{\rm dim}\,H_{2s+1}(M,k)$
and let ${\mathfrak m }:=H_{2s}(M,k)\oplus H_{2s+1}(M,k)\oplus
H_{4s+1}(M,k)$ be the maximal ideal of $R$. Given a finitely
generated $R$-module $M$, let $M_{0}:=M/{\mathfrak m}M$. $M_{0}$
is a finite dimensional vector space over $k$. Below we will use
the following form of Nakayama's lemma: if $\varphi:M\to N$ is a
morphism of $R$-modules such that the induced morphism
$\varphi^{0}:M_{0}\to N_{0}$ is surjective, then $\varphi$ is also
surjective.

To compute $\mbox{\rm Tor}^{R}(k,k)$ we need to take a projective
resolution of $k$ regarded as a $R$-module in the obvious way.
Since $R$ is local a $R$-module is projective if and only if is
free. Hence, we will construct a resolution of the form:

\[\cdots\rightarrow R^{b_{i}}{\buildrel\partial_{i}\over\rightarrow}\cdots\rightarrow R^{b_{1}}
{\buildrel\partial_{1}\over\rightarrow} R{\buildrel\partial_{0}
\over \rightarrow}k\rightarrow 0.\]

The first map $\partial_0$ is given simply by
\[\partial_{0}(x,y,z,t)=x,\]
where $(x,y,z,t)\in R=H_{0}(M,k)\oplus H_{2s}(M,k)\oplus
H_{2s+1}(M,k)\oplus H_{4s+1}(M,k)$ and we identify $H_{0}(M,k)$
with $k$. Clearly ${\rm Ker}\,\partial_{0}={\mathfrak m}$.

We will now define a surjective morphism $\partial_{1}:R^{2a}\to
{\mathfrak m}$. Let ${\bf 1}=(1,0,0,0)\in R$. Clearly ${\bf 1}$
generates $R$ and hence given any free module $R^{b}$, the
elements ${\bf e}_{i}=(0,\dots,{\bf 1}_{i},\dots,0)$
 for $1\leq i\leq b$ generate $R^{b}$. Hence, to define $\partial_{1}$ it suffices to indicate
the images of the ${\bf e}_{i}$'s. Pick a basis of
$H_{2s}(M,k)\oplus H_{2s+1}(M,k)$ (which has dimension $2a$) and
let $\partial_{1}$ be determined by a bijection between the
generators of $R^{2a}$ and this basis.

Note that ${\mathfrak m}_{0}={\mathfrak m}/{\mathfrak m}^{2}\cong
H_{2s}(M,k)\oplus H_{2s+1}(M,k)$. Hence $\partial_{1}^{0}$ is an
isomorphism and by Nakayama's lemma $\partial_{1}$ is surjective.

Let $Q\subset R$ be the ideal given by those elements of the form
$(0,0,0,t)$. Note that
\begin{enumerate}
\item ${\mathfrak m}{\rm Ker}\,\partial_{1}=Q^{2a}$;
\item  ${\rm Ker}\,\partial_{1}/{\mathfrak m}{\rm Ker}\,\partial_{1}$ has dimension $4a^{2}-1$.
\end{enumerate}

To define $\partial_{2}$, we take $R^{4a^{2}-1}$ and we map the
canonical $4a^{2}-1$ generators of $R^{4a^{2}-1}$ onto a basis of
${\rm Ker}\,\partial_{1}/{\mathfrak m}{\rm Ker}\,\partial_{1}$.
This gives a surjective morphism as before.

By continuing in this fashion we find that at the i-th step of the
construction of the resolution we have:

\begin{enumerate}
\item ${\mathfrak m}{\rm Ker}\,\partial_{i-1}=Q^{b_{i-1}}$;
\item  ${\rm Ker}\,\partial_{i-1}/{\mathfrak m}{\rm Ker}\,\partial_{i-1}$ has dimension $2ab_{i-1}-b_{i-2}$.
\end{enumerate}

Therefore $b_{i}=2ab_{i-1}-b_{i-2}$. This implies that the growth
of sequence $b_{i}$ is exponential if $a>1$ (with exponent
$a+\sqrt{a^{2}-1}$) and at most linear if $a\leq 1$.

Now observe that we have the isomorphism
$R^{b_{i}}\otimes_{R}k\cong k^{b_{i}}$ and under this isomorphism
the map $\partial_{i}\otimes 1$ is zero. Thus the differential of
the complex $R^{b_{i}}\otimes_{R}k$ is zero, so the dimensions of
 $\mbox {\rm Tor}^{R}(k,k)$ over $k$ grow exactly as the
$b_{i}$'s.

\end{proof}

Closed simply connected smooth 5-manifolds have been classified by
S. Smale in the spin case \cite{smale} and by D. Barden
\cite{Barden} in the general case. We will now briefly describe
the classification.

The oriented (5-dimensional) cobordism group has order 2. The
non-trivial cobordism class is formed by the manifolds for which
the Stiefel-Whitney number $w^2 \cup w^3 \neq 0$. Let $M$ be a
closed simply connected smooth 5-manifold. If $M$ bounds, then the
torsion part of $H_2 (M,\Z)$ is isomorphic to $G\oplus G$ for some
finite Abelian group $G$. If $M$ belongs to the non-trivial
cobordism class then the torsion part of its  second homology
group is of the form ${\Z}_2 \oplus G \oplus G$, where $G$ is
again a finite Abelian group.

The second Stiefel-Whitney class of a simply connected closed
manifold is given by a homomorphism $w^2 :H_2 (M, \Z )\rightarrow
{\Z}_2$. There exists a basis of the Abelian group $H_2 (M, \Z )$
such that it has the maximal possible number of elements (for a
basis of the Abelian group) and such that $w^2$ does not vanish in
at most one of the elements of the basis. If the order of this
element is $2^{i}$ then $i$ depends only on $M$.

This invariant $i(M)$ together with $H_2 (M,\Z )$ is a complete
set of invariants for simply connected closed 5-manifolds.

\vspace{.3cm}

Let $X_{-1} =SU(3)/SO(3)$ be the Wu-manifold, which is
characterized by $i(X_{-1})=1$ and $H_2 (X_{-1} ,\Z )= {\Z }_2$.
Let $X_0 =S^5$, $M_{\infty} = S^3 \times S^2$ and $X_{\infty}
={\eta}_3$ (the only non-trivial $S^3$-bundle over $S^2$).

For $1\leq j<\infty$ let $X_j$ be a closed simply connected
non-spin 5-manifold such that $H_2 (X_j ,\Z )={\Z}_{2^j} \oplus
{\Z}_{2^j}$. Then $i(X_j )=j$. Also let $M_j$ be a spin manifold
with $H_2 (M_j ,\Z )={\Z}_j \oplus {\Z}_j$. Of course, $i(M_j)=0$.

Then Barden proves that any simply connected closed 5-manifold $M$
is diffeomorphic to a connected sum of some of these manifolds.
More precisely, $M=X_j \# M_{k_1} \# ...\#M_{k_l}$ where $-1\leq
j\leq \infty $, $k_1 >1$ and $k_i$ divides $k_{i+1}$ for all $i$.
Note that then $i(M)=j$ and $H_2 (M,\Z )= {\Z}_{2^j} \oplus
{\Z}_{2^j} \oplus {\Z}_{k_1} \oplus {\Z}_{k_1} \oplus ... \oplus
{\Z}_{k_s} \oplus {\Z}_{k_s}$, unless $j=-1$ in which case the
first two factors should be replaced by one copy of ${\Z}_2$.

As a consequence of Theorem \ref{loop5m} and the classification of
simply connected 5-manifolds we obtain:

\newpage

\begin{Corollary} Let $M$ be a closed
simply connected 5-manifold. Then $M$ is elliptic if and only if $M$ is diffeomorphic
to:

\begin{enumerate}
\item $S^{5}$;
\item $S^{3}\times S^{2}$ whose second homology group is $\Z$ and is a spin
manifold;
\item $\eta_{3}$, the nontrivial $S^{3}$-bundle over $S^{2}$, whose second homology group is $\Z$ and is not spin;
\item the Wu-manifold $X_{-1}=SU(3)/SO(3)$ whose second homology group is $\Z_{2}$ and is not spin.
\end{enumerate}
Moreover if $M$ is not elliptic, the homology of the loop space of
$M$ grows exponentially for some field of coefficients $k_{p}$.
\label{class5m}
\end{Corollary}

\begin{table}
\caption{The elliptic list in dimensions 4 and 5}

\begin{tabular}{|c|c|}
\hline {\rm dim} 4  &  {\rm dim} 5\\
\hline
          &       \\
$S^{4}$   & $S^{5}$\\
&\\
 $\cp^{2}$  & $S^{3}\times S^{2}$\\
&\\
 $S^{2}\times S^{2}$ &  $X_{-1}=SU(3)/SO(3)$\\
&\\
 $\eta_{2}=\cp^{2}\#\overline{\cp}^{2}$ & $\eta_{3}$\\
&\\
 $\cp^{2}\# \cp^{2}$ &  \\
\hline
\end{tabular}

\end{table}

\section{Existence of a metric with zero entropy on each manifold in the elliptic list}
\label{zeroentropy}
\subsection{Dimension 4}

The standard symmetric metrics on $S^{4}$ and $\cp^{2}$ have all
the geodesics closed and with the same period, and hence
their geodesic flows have zero topological entropy.
On $S^{2}\times S^{2}$ consider the product metric of the round metric
on $S^{2}$; it follows from part (1) in Lemma \ref{facil}
that the geodesic flow of the product metric has zero entropy.

The manifold $\cp^{2}\#\overline{\cp}^{2}$ is $\eta_{2}$, the non-trivial
$S^{2}$-bundle over $S^{2}$, and it is known to be diffeomorphic
to the space that we now describe.
Represent $S^{3}\subset \C^{2}$ as
pairs of complex numbers $(z_{1},z_{2})$ with $|z_{1}|^{2}+|z_{2}|^{2}=1$.
Let $S^{1}$ act on $S^{3}$ by
\[(w,(z_{1},z_{2}))\mapsto (wz_{1},wz_{2}),\]
where $w\in S^{1}$ is a complex number with modulus one.
Let $S^{1}$ also act on $S^{2}$ by rotations.
Consider the space $M=S^{3}\times_{S^{1}}S^{2}$
obtained by taking the quotient of $S^{3}\times S^{2}$ by the diagonal
action of $S^{1}$. The manifold $M$ is diffeomorphic
to $\cp^{2}\#\overline{\cp}^{2}$.
Endow $S^{3}$ and $S^{2}$
 with the canonical metrics of curvature one.
By part (1) of Lemma \ref{facil}
the product metric on $S^{3}\times S^{2}$ has zero entropy.
By part (2) in Lemma \ref{facil} the submersion metric
on $M=S^{3}\times_{S^{1}}S^{2}$ will also have a geodesic flow
with zero entropy.

We are left with the case of $M=\cp^{2}\#\cp^{2}$ which is in fact the only tricky case.
The manifold $M$ can be obtained from two copies
of $S^{3}\times_{S^{1}}D^{2}$ where $D^{2}$ is the 2-disk
and $S^{1}$ acts diagonally, glued along their boundary
$S^{3}\times_{S^{1}}S^{1}=S^{3}$ by an orientation reversing map.
In \cite{P1} the first author proved
that the metrics considered by J. Cheeger in \cite{Ch} have zero topological entropy.

\subsection{Dimension 5} The round metric on $S^{5}$ and the product metric
on $S^{3}\times S^{2}$ clearly have zero entropy.

For the Wu manifold $X_{-1}$ we proceed as follows. Let us consider a biinvariant metric on
$SU(3)$. Since every geodesic is the orbit of a 1-parameter subgroup and since $SU(3)$ is compact
it follows easily that all the Jacobi fields grow at most linearly.
Therefore all the Liapunov exponents of {\it every} geodesic in $SU(3)$ are zero. It
follows from Ruelle's inequality \cite{Ru} that all measure
entropies are zero.  Hence, by the variational principle, the
topological entropy of the geodesic flow of $SU(3)$ must be zero. Endow $X_{-1}$ with the submersion
metric. It follows from part (2) in Lemma \ref{facil} that this metric has zero
topological entropy.

We are left with $\eta_{3}$. This is handled in a similar way with the help of the next
lemma which gives a convenient way of expressing $\eta_{3}$ using group actions.

\begin{Lemma}Consider on $S^{3}\times S^{3}\subset \C^{2}\times C^{2}$ the action of
$S^{1}$ given by
\[(w,(z_{1},z_{2},z_{3},z_{4}))\mapsto (wz_{1},wz_{2},wz_{3},z_{4}),\]
where $w\in S^{1}$ is a complex number with modulus one. This action is fixed point
free and the quotient of $S^{3}\times S^{3}$ by this action is $\eta_{3}$.
\end{Lemma}

\begin{proof}Let $M$ be the quotient of $S^{3}\times S^{3}$ by the circle action.
A simple argument with the long exact sequence of the fibration
shows that $M$ is simply connected and $\pi_{2}(M)=\Z$. By the
Hurewicz theorem $H_{2}(M,\Z)=\Z$. Note that $M$ contains a copy
of $\cp^{2}\#\overline{\cp}^{2}$ given by the projection to $M$ of
the subset of $S^{3}\times S^{3}$ given by $\{\mbox{\rm imaginary
part of }z_{4}=0\}$ and hence $M$ is not spin. It follows from the
Barden-Smale classification that the only closed simply connected
non spin 5-manifold with $H_{2}(M,\Z)=\Z$ is $\eta_{3}$.

\end{proof}

Using the lemma it is easy to construct a metric on $\eta_{3}$ with zero entropy.
Consider on $S^{3}\times S^{3}$ the product metric and on the quotient the submersion metric.
By Lemma \ref{facil} the metric thus constructed on $\eta_{3}$ has zero entropy.

\section{$\mathcal{F}$-structures and minimal entropy}
\label{F-str}

We begin by considering the case of a non-trivial $S^1$-action.
This preliminary result will not be used in the proof for the
case of a general $\mathcal{F}$-structure. But we think that its
much simpler proof gives a nice picture of the ideas behind the
general case.

\begin{Theorem} Suppose that the closed connected smooth
manifold $M$ admits a non-trivial $S^1$-action. Then the
minimal entropy of $M$ is 0.
\end{Theorem}

\begin{proof} First consider a metric $g$ on $M$ which is
invariant under the $S^1$-action. This is obtained as usual
by averaging any given Riemannian metric over the orbits.

Now consider the manifold $\bar{M} = M\times S^1$ and
for any $\delta >0$ the Riemannian
metric ${\bar{g}}_{\delta} =g+ \delta dt^2$ (where
$dt^2$ is the Euclidean metric on $S^1$) on $\bar{M}$.

Define a (free)  $S^1$-action on $\bar{M}$ by

$$\lambda . (x,\theta )= (\lambda ._M x, \lambda \theta)$$

The quotient of $\bar{M}$ by this action is diffeomorphic to $M$
and the metric ${\bar{g}}_{\delta}$ is invariant through the action;
therefore
it induces a metric $g_{\delta}$ on $M$. The  projection

$${\pi}:(\bar{M} , {\bar{g}}_{\delta} ) \rightarrow (M,g_{\delta} )$$

\noindent
is a Riemannian submersion,
and therefore the entropy of $g_{\delta}$ is bounded above by
the entropy of ${\bar{g}}_{\delta}$
(see Section 2.5) which is actually equal to
the entropy of $g$
(see Section 2.5). Therefore to prove the theorem it is enough
to show that the volume of $(M,g_{\delta})$ approaches 0 as
$\delta$ approaches  0. We will prove this now.

First we identify the quotient
(of $\bar{M}$ by the $S^1$-action) with $M$ via the diffeomorphism
which sends $x\in M$ to the class of $(x,1)$. Let $v_x$ be the
vector tangent to the $._M$-action at $x$ and let $\omega$ be
the tangent to the canonical $S^1$-action on $S^1$ (which gives
the usual trivialization of the tangent space of $S^1$).
Let ${\varepsilon}_x =g(v_x ,v_x )$.

The tangent vector to the action on $\bar{M}$ is $(v_x ,\omega )$.
If ${\varepsilon}_x \neq 0$,
the ${\bar{g}}_{\delta}$-orthogonal subspace to this vector is
spanned by $(v_x ,-\frac{{\varepsilon}_x}{\delta } \omega )$ and the
subspace of vectors of the form $(v,0)$ where
$v\in V_x \subset T_x M$, the subspace of vectors
$g$-orthogonal to $v_x$. It is clear that
$g_{\delta}$ and $g$ coincide on $V_x$.
Moreover, $g_{\delta} (v_x ,v)=0$ for all $v\in V_x$.

Since

$$(v_x ,0)= \frac{{\varepsilon}_x /\delta}{1+({\varepsilon}_x /\delta )}
\left( v_x ,\omega \right)
+   \frac{1}{1+({\varepsilon}_x /\delta )} \left( v_x
  ,-({\varepsilon}_x /\delta )\omega \right),$$

\noindent
we have that

$$g_{\delta} (v_x ,v_x )=
{\left( \frac{1}{1+({\varepsilon}_x /\delta ) }\right)}^2
  \ {\bar{g}}_{\delta}
\left( (v_x  ,-({\varepsilon}_x /\delta ) \omega ),
                   (v_x ,-({\varepsilon}_x /\delta ) \omega ) \right) $$

$$=\frac{{\varepsilon}_x +
{\varepsilon}_x^2 /\delta }{(1+({\varepsilon}_x  /\delta ))^2}=
\frac{\delta}{\delta +{\varepsilon}_x} \ \ g(v_x ,v_x ).$$

This implies  the following equation for the volume elements
of the two metrics:

$$\mbox{\rm dvol}(g_{\delta})=
  \sqrt{\delta}\frac{1}{\sqrt{\delta +\varepsilon}}\mbox{\rm dvol}(g).$$

This formula will be enough to show that the volume of
$(M, g_{\delta})$ approaches 0 with $\delta$.
Note first that the formula shows
that the volume of any region computed with $g_{\delta}$ is
always at most the
volume of the same region computed with $g$ (independently of $\delta$).
Given any $\rho >0$, we can find an open neighborhood of the fixed
point set of the $S^1$-action on $M$ which has $g$-volume less
than $\rho /2$. Then the $g_{\delta}$-volume of this neighborhood
will also be less than $\rho /2$ for any $\delta$. Away from the
neighborhood, $\varepsilon$ has a positive lower bound, and the
volume formula clearly shows that the $g_{\delta}$-volume of
the complement of the neighborhood is of the order of
$\sqrt{\delta}$ for $\delta$ small.
Therefore, for $\delta$ small enough the volume of the complement
will also be less than $\rho /2$.
This completes the proof
of the theorem.

\end{proof}

This result should be compared to the collapsing with bounded
sectional curvature of Cheeger and Gromov \cite{Cheeger, Gromov}.
If the manifold $M$ admits a
locally free $S^1$-action then picking a Riemannian metric
$g$ on $M$ invariant through the action and then shrinking
along the orbit produces a sequence of metrics with uniformly
bounded curvature and volume and injectivity radius converging
to zero. This is not true if the action has fixed points. To
get a geometrical picture of our theorem one can consider the
canonical $S^1$-action on $S^2$ which has fixed points in the
poles. The metrics produced in the proof of the theorem
will shrink the horizontal circles by a non-constant factor,
which approaches 1 near the poles.
At the poles the curvature will blow-up and the injectivity radius
will stay uniformly bounded from below. But the volume will
collapse and the entropy will stay bounded.

Cheeger and Gromov introduced in \cite{Cheeger, Gromov} the
concept of ${\mathcal F}$-structures and generalize the previous
result to manifolds admitting ${\mathcal F}$-structures with
certain special properties: polarized $\mathcal{F}$-structures of
positive rank. There exist plenty of examples of closed manifolds
which admit $\mathcal{F}$-structures but which cannot be collapsed
with bounded sectional curvature; manifolds whose minimal volume
is non-zero. Therefore they do not admit polarized
$\mathcal{F}$-structures of positive rank. We will show that the
minimal entropy does vanish in the presence of general ${\mathcal
F}$-structures. We will follow the notation of \cite{Cheeger} as
closely as possible, and the reader should check that reference
for any detail about the definition and many constructions related
to $\mathcal{F}$-structures we will use. We consider first the
definition of an $\mathcal{F}$-structure.

A sheaf of tori $\mathcal{S}$ over the smooth
manifold $M$ is said to {\it act} on $M$ if for each open subset
$U$ of $M$ there is a local action of
the group of sections $\mathcal{S} (U)$ on $U$,
with the obvious compatibility between restriction
homomorphisms of the sheaf and restrictions of the local actions
(a local action
of a group $G$ is an action defined only on a neighborhood of
$\{ e \} \times U \subset G\times U $). The action divides $M$ into
{\it orbits} and a subset of $M$ is called {\it saturated} if it is
a union of orbits.

\begin{Definition}
An $\mathcal{F}$-structure on a smooth closed manifold $M$ is
given by an action on $M$ of a sheaf $\mathcal{S}$ of tori
together with a finite cover of $M$ by saturated open subsets $\{
U_1 ,...,U_N \}$ such that:

\noindent
{\bf (a)} On each $U_i$ there is a locally constant subsheaf
${\mathcal{S}}_i$ of $\mathcal{S}$ and a finite normal covering
${\pi}_i  : \widetilde{U_i}   \rightarrow  U_i$ such that
the structure homomorphisms
of ${{\pi}_i}^* ({\mathcal{S}}_i )$
give  isomorphisms between the global sections and the stalks.

\noindent
{\bf (b)} The local action of the sections defines a smooth,
effective torus action

$$._i :T^{k_i} \times \widetilde{U_i} \rightarrow \widetilde{U_i} ,$$

\noindent
{\bf (c)} The stalk of the sheaf at any $x\in M$ is spanned by the stalks
of the subsheaves corresponding to the $U_i$'s which contain
$x$ and
non-empty intersections of the $U_i$'s
also have a finite covering such that the pull back of the sheaf
spanned by the corresponding ${\mathcal{S}}_i$'s
gives rise to a global torus action as before.
\end{Definition}

\begin{Definition} An $\mathcal{F}$-structure is called a
$\T$-structure if all the coverings ${\pi}_i :\widetilde{U_i}
\rightarrow U_i$ are trivial.
\end{Definition}

\begin{Remark}{\rm The dimension of the orbit through $x$ 
is called the {\it rank} of
$\mathcal{F}$ at $x$. The minimum of the dimensions of the orbits
is called the {\it rank} of the $\mathcal{F}$-structure. The
$\mathcal{F}$-structure is called {\it polarized} if the torus
actions defined on the finite coverings are locally-free . }
\end{Remark}

\begin{Remark}{\rm Our definition of $\mathcal{F}$-structure is
essentially the same as the one in \cite{Cheeger}. More precisely,
one can see that given any $\mathcal{F}$-structure as defined
by Cheeger and Gromov there exists an {\it atlas} with the properties
in our definition (see page 317 in \cite{Cheeger}). }
\end{Remark}

\begin{Remark} {\rm A $\T$-structure
is given by a covering by open subsets
and a torus action on each subset such that any
intersection of the open subsets is invariant through the corresponding
actions and these commute. The stalk over any point $x$
of the sheaf appearing in Definition 5.2 is the maximal
torus which is acting on $x$. 
The definition is of course the same as
the original one given by Gromov in \cite{Gromov}, except that it is
only asked  that the torus actions are effective (but not necessarily
locally free). }
\end{Remark}

\begin{Example}{\rm Any non-trivial $S^1$-action on $M$ is of course a
$\T$-structure on $M$. Hence, for instance, $S^4$ and $\cp^2$
admit $\T$-structures although they cannot admit any polarized
$\mathcal{F}$-structure.}
\end{Example}

\begin{Example}{\rm The compact complex surface $K3$ admits a $\T$-structure, even
though it does not admit any non-trivial $S^1$-action
\cite{AH}. Actually every elliptic compact complex surface
admits a $\T$-structure as we will show below.}
\end{Example}

We will see now that $\T$-structures behave relatively well with
respect to the usual operations of connected sums and surgeries on
manifolds. T. Soma proved in \cite{Soma} that the family of
3-manifolds which admit polarized $\T$-structures is closed under
connected sums. As pointed out  by Gromov in \cite{Gromov}, this
result generalizes to any odd dimension. We will see now that the
result also holds for the family of manifolds which admit general
$\T$-structures and for any dimension greater than 2.

\begin{Theorem} Suppose $X$ and $Y$ are $n$-dimensional
manifolds, $n>2$, which admit a $\T$-structure. Then $X\# Y$ also
admits a $\T$-structure.\label{connectedsums}
\end{Theorem}

\begin{proof} Pick a point $x\in X$ so that $x$ lies in only one of
the open subsets of the $\T$-structure (for this one might need to
do some harmless changes in the $\T$-structure, like eliminating
any open subset which is contained in the union of the others). We
can also assume that the torus acting on the open subset
containing $x$ is of dimension one and that $x$ lies on a regular
orbit.

Now pick a small $(n-1)$-ball $D_x$ centered at $x$ and transverse to the
$S^1$-action.
The union of the orbits through $D_x$ form an embedded solid torus
$S^1 \times D_x$. Repeat the same procedure to obtain an embedded
solid torus $S^1 \times D_y$ in $Y$ containing a point $y\in Y$.
We will perform the connected sum inside $S^1 \times D_x$ and
$S^1 \times D_y$.

First divide $D_x$ into an inner ball and an outer annulus: $D_x
=D_{{\varepsilon}_1} \cup (S^{n-2} \times [{\varepsilon}_1
,{\varepsilon}_2 ])$. We can identify $S^1 \times D_x \# S^1
\times D_y$ with $S^1 \times D_x -S^{n-2} \times D^2$, where $D^2$
is a small 2-dimensional ball centered at a point in the middle of
$S^{n-2} \times [{\varepsilon}_1 ,{\varepsilon}_2 ]$ and
transverse to $S^{n-2}$ in $S^1 \times S^{n-2} \times
[{\varepsilon}_1 , {\varepsilon}_2 ]$. The component of the
boundary corresponding to the boundary of the deleted $S^{n-2}
\times D^2$ is identified with the boundary of $S^1 \times D_y$.


We can now describe the $\T$-structure on $X \# Y$. On $(X-S^1
\times D_x )\cup S^1 \times D_{{\varepsilon}_1}$ leave the initial
$\T$-structure. On $S^1 \times S^{n-2} \times [{\varepsilon}_1
,{\varepsilon}_2 ] - (S^{n-2} \times D^2 )$ consider any
non-trivial $S^1$-action  on the $S^{n-2}$-factor (here is where
we need the hypothesis $n>2$). The action induced on each
component of the boundary glues to the canonical action on the
$S^1$-factor to create a $T^2$-action (in case $n$ is even it will
have orbits of dimension 1). Finally on $Y-(S^1 \times D_y )$
leave the initial $\T$-structure.

\end{proof}

\begin{Theorem}
Every compact complex elliptic surface admits a $\T$-structure.
\label{telliptic}
\end{Theorem}

\begin{proof} For the proof we will need smooth descriptions of 
the surfaces: see
\cite{Friedman, LeBrun2} for details. Every elliptic surface of
Euler characteristic 0 is obtained by performing logarithmic
transforms on a {\it basic} elliptic surface. Every elliptic
surface is obtained by taking the {\it fiber sum} of an elliptic
surface of Euler characteristic 0 and some rational elliptic
surfaces, and then blowing up some points.

Basic surfaces are fiber bundles with fibers $T^2$ and structure
group in $SL(2,\Z)$. Hence they admit a polarized $\T$-structure
whose orbits are the fibers.

Now let $B\times T^2$ be a neighborhood of a fiber on a basic
surface $M$, where $B$ is identified with the unit ball in
$\C={\re}^2$. Fix a positive integer $m$ and integers $a$, $b$
such that $(a,b)$ has order $m$ in $\Z_m \oplus \Z_m$. Let
$F:B\times T^2 \rightarrow B\times T^2$ be given by
$F(z,t)=(e^{2\pi i/m}z,t)$. $F$ generates a group $G_1$ of
diffeomorphisms of $B\times T^2$ of order $m$. The quotient of
$B\times T^2$ by this group is again diffeomorphic to $B\times
T^2$. Consider also the map $L:B\times T^2 \rightarrow B\times
T^2$ given by $$L(z,t) = (e^{2\pi i /m}z, (t_1 e^{2a\pi i/m},t_2
e^{2b\pi i/m})).$$

\noindent
$L$ generates a group $G_2$ of
diffeomorphism of $B\times T^2$ of order $m$ which acts freely on
$B\times T^2$. The map
$P:S^1 \times T^2 /G_1  \rightarrow S^1 \times T^2 /G_2$,

$$P(z,t)= (z, (z^a t_1 ,z^b t_2 ))$$

\noindent is a diffeomorphism. The logarithmic transform (of order
$m$) at the fiber over $(0,0)$ in $M$ is the 
elliptic surface $\widetilde{M}$
obtained by gluing $M-B\times T^2$ and $B\times T^2 /G_2$ via this
diffeomorphism. Clearly the obvious $S^1$-action on $B$ (which
fixes $(0,0)$) induces an $S^1$-action on $B\times T^2 /G_2$ which
commutes with the action on the fibers. Hence $\widetilde{M}$
admits a $\T$-structure (with orbits of dimension 0,1,2 and 3).

Rational elliptic surfaces are diffeomorphic to $S=\cp^2 \#
9{\overline\cp}^2$ and therefore admit $\T$-structures by the
previous theorem. Nevertheless we will need to perform fiber sums
and so we will give another $\T$-structure on it, compatible with
the elliptic fibration. To do this we need first to give a
description of the surface as an elliptic surface (see
\cite{LeBrun2}). Let $T^2 ={\re}^2 /\Z^2$ and consider the
involution $I(z)=-z$ of $T^2$. Let $H:S^2 \rightarrow S^2$ be
rotation of $180^{\circ}$ around the $z$-axis. The diffeomorphism

$$J=(I,H):T^2 \times S^2 \rightarrow T^2 \times S^2$$

\noindent has 8 fixed points. Identify a neighborhood of each of
these points with a ball $B$ in $\C^2$. Consider $U=\{ (z,l) \in
B\times \cp^1 :  \  z\in l \}$. The canonical projection
${\pi}_1 |_U : U\rightarrow B$ induces an isomorphism away from
the preimage of $0$. Construct a surface $\widetilde{S}$ by
replacing the eight copies of $B$ with $U$ in $T^2 \times S^2$.
The involution $J$ extends to an involution $\widetilde{J}$ on
$\widetilde{S}$ which has 8 spheres as the set of fixed points.
Then $S=\widetilde{S} /\widetilde{J} $. Let $\pi :S^2 \rightarrow
S^2 =S^2 /H$ be the projection. Then $\pi \circ {\pi }_2 :T^2
\times S^2 \rightarrow S^2$ induces a map 
$p:\widetilde{S} \rightarrow S^2$ which commutes with $\widetilde{J}$
and so induces a map 
$S\rightarrow S^2$ whose
generic fiber is $T^2$; this map expresses $S$ as an elliptic
surface.

Note that the $S^1$-action  on $B$ given by 
$\lambda (w_1 ,w_2 )= (w_1 ,\lambda w_2 )$
commutes with $J$ and induces an $S^1$-action on $U$.
We can extend this action to an $S^1$-action defined 
on a neighborhood of the fibers of $p$ over the north 
and south poles. This actions commutes with $\widetilde{J}$
and so induces an action on a neighborhood of the fibers
of $S\rightarrow S^2$ over the poles.

Away from the fibers over the poles $S$ actually is the total
space of a fiber bundle with structure group $\{ Id, I\}$. There
is then a polarized  $\T$-structure defined on this piece, whose
orbits are the fibers. 
On the boundary of the neighborhoods around the fibers over
the poles the two actions commute.
This defines a $\T$-structure on $S$.

The fiber sum of two elliptic surfaces is done as follows: pick
regular fibers on each surface identifying neighborhoods of them
with $D\times T^2$ ($D$ is a small 2-ball). Delete the
corresponding regular fiber from each surface and then glue both
surfaces along $(D\# D)\times T^2$. The diffeomorphism class of
the resulting surface will depend only on the isotopy class of the
diffeomorphism chosen to identify the fibers with $T^2$. We can
therefore take the diffeomorphism to be in $SL(2,\Z)$ and we can
see that the $\T$-structures we defined on the surfaces of Euler
characteristic 0 and the rational elliptic surfaces glue well
along the fiber sum.

Finally blowing up points means, in terms of diffeomorphisms, to
take connected sums with ${\overline\cp}^2$'s. Such a connected
sum admits a $\T$-structure by the previous theorem.

\end{proof}

We can now also see that inside the family of manifolds with
$\T$-structures one can perform surgery on spheres which are
``well positioned'' with respect to the $\T$-structure.

\begin{Definition}Let $M$ be a manifold with a fixed $\T$-structure.
An embedded $k$-sphere ${\bf{S}}^k$ is said to be {\it completely
transversal} with respect to the $\T$-structure if:

\noindent
{\rm 1)} ${\bf{S}}^k$ intersects only one of the open subsets of the
$\T$-structure.

\noindent
{\rm 2)} The torus acting on the open subset of (1) has dimension 1 and
the orbits passing through ${\bf{S}}^k$ form an embedded $S^k
\times S^1$ with trivial normal bundle.
\end{Definition}

\begin{Remark}{\rm Note in particular that the normal bundle of a completely
transversal sphere is trivial.}
\end{Remark}

\begin{Example}{\rm If $X$ admits a $\T$-structure and $Y$ is any other manifold
then $X\times Y$ admits an obvious $\T$-structure. Any homotopy
class in $Y$ which can be represented by an embedded sphere with
trivial normal bundle (in $Y$) can be represented by a completely
transversal sphere (in $X\times Y$).}
\end{Example}

\begin{Theorem} Let $M^n$ be a manifold with a $\T$-structure.
Let ${\bf{S}}^k$ be a completely transversal sphere (with respect
to the given $\T$-structure). The manifold $\hat{M}$, obtained by
performing surgery on $\bf{S}$, also admits a $\T$-structure.
Moreover, if $n$ and $n-k$ are odd and the structure on $M^n$ is
polarized, then $\hat{M}$ also admits a polarized $\T$-structure.
\label{surgery}
\end{Theorem}

\begin{proof}
Let $S^k \times S^1 \times D^{n-k-1}$ be a tubular neighborhood of
the union of the orbits through $\bf{S}$. Consider the unit
$n$-sphere $S^n \subset {\re}^{n+1}$. Pick a non-trivial
$S^1$-action on $S^n$, for instance complex multiplication in the
first 2 coordinates. In case $n$ is odd we can pick a free
$S^1$-action. Choose a regular orbit of the action and a disc
$D^{n-1}$ transverse to the orbit. Pick a canonical embedded
$k$-sphere $S^k_0 \subset D^{n-1}$ and a tubular neighborhood
$S^k_0 \times S^1 \times D^{n-k-1}$ of the union of the orbits
through $S^k_0$. The manifold $\hat{M}$ is obtained by gluing $M$
and $S^n$ along $\bf{S}$ and $S^k_0$. But gluing two copies of
$S^k \times S^1 \times D^{n-k-1}$ along the $k$-spheres is the
same as taking the product of a $k$-sphere with the connected sum
of two copies of $S^1 \times D^{n-k-1}$. Hence in this glued part
we can consider the  $\T$-structure we defined in the previous
theorem, which on each component of the boundary coincides with
the structure of $M$ and $S^n$, respectively. This clearly defines
a $\T$-structure on $\hat{M}$. This structure is polarized if $n$
and $n-k$ are odd.

\end{proof}

\section{Collapsing with bounded entropy: Proof of Theorem A}
\label{zerome}

In this section we will prove that the minimal entropy of a
closed manifold which admits an
$\mathcal{F}$-structure vanishes. The general idea of the proof
is quite simple. Given an $\mathcal{F}$-structure on $M$ we
define a polarized $\mathcal{F}$-structure
on $M\times T^k$ for
some $k$ and consider a Riemannian metric on the product which
is invariant through all the torus actions.
Then we collapse the metric along the orbits of
the $\mathcal{F}$-structure on $M\times T^k$. The procedure
constructs
metrics which are invariant by the canonical $T^k$-action on
$M\times T^k$. Taking the quotient by this action
gives a Riemannian submersion over a metric
on $M$. Now, for the polarized structure on $M\times T^k$, Cheeger
and Gromov \cite{Cheeger} proved that the sectional curvatures of
the collapsed metrics are uniformly bounded. Therefore the entropy
of the metrics are also uniformly bounded (see Section 2.3)
and since entropy is
non-increasing under Riemannian submersions
(see Section 2.5), the collapsed metrics
on $M$ also have uniformly bounded entropy. The theorem therefore
reduces to the proof that the volumes of the metrics on $M$ collapse.
Note that the only properties about entropy we will use in the proof
are its bound in terms of curvature and
its behaviour under Riemannian submersions. Since Riemannian
submersions do not decrease sectional curvatures, the same
proof works for any quantity that depends only on lower bounds
for the sectional curvature. We will use this remark in the
next section to study certain curvature invariants for manifolds
admitting $\mathcal{F}$-structures.

In the proof of the theorem we will need the following
elementary lemma from linear algebra:

\begin{Lemma} Let $(V_1 ,h_1 )$ and $(V_2 ,h_2)$ be two real vector
spaces of dimension $l$ with inner products. Let $F$ be a subspace
of $V_1 \oplus V_2$ of dimension $l$ which intersects trivially
with  both
$V_1$ and $V_2$ such that for any $(v,w)\in F$,
$h_1 (v,v)\leq h_2 (w,w)$. Then:

\noindent
{\bf a.} Consider $F$ as the graph of a map $\widetilde{F} :V_2 \rightarrow
V_1$ and let $I :V_1 \rightarrow V_1$  be given by  $I (v)={\pi}_1 \circ
{\pi}_F (v,0)$ (${\pi}_F :V_1 \oplus V_2 \rightarrow F$ is the orthogonal
projection). Then $(\det I )^2 \geq 4^{-l} (\det \widetilde{F} )^{4l}$.

\noindent
{\bf b.} Given any $\lambda$, $0<\lambda \leq 1$ consider
the inner product ${\bar{h}}_{\lambda}$ on $V_1 \oplus V_2$ defined
by $\lambda (h_1 +h_2 )|_{F} +(h_1 +h_2 )|_{F^{\perp}}$. Let
$h_{\lambda}$ be the inner product on $V_1$ obtained as the
quotient of ${\bar{h}}_{\lambda}$ (by ${\pi}_1$).
Then $dvol(h_{\lambda} )
\leq  dvol(h_1 )$.

\end{Lemma}

\begin{proof} {\bf a)} Consider an orthonormal basis
$\{ v_1 ,...,v_l \}$ of $(V_1 ,h_1 )$.
If $I (v_j )=a_{ij} v_i $ and we let $A=(a_{ij} )$ then
$(\det I )^2 = \det (A^t A)$. But $A^t A$ is a positive definite
symmetric matrix, and therefore it has $l$ positive eigenvalues
${\mu}_1 \leq {\mu}_2 \leq ...\leq {\mu}_l $ and
$(\det I )^2  ={\mu}_1 ...{\mu}_l \geq {{\mu}_1 }^l$.

Now consider an orthonormal basis $\{ w_1 ,...,w_l \}$ of
$(V_2 ,h_2 )$ and let $\widetilde{F} (w_j )=b_{ij} v_i$.
Let $\gamma =(\det \widetilde{F} )^2$.
If $B=(b_{ij} )$ then $\gamma = \det (B^t B)$. Again,
$B^t B$ is a positive definite symmetric matrix. Moreover,
since $h_1 (v,v) \leq h_2 (w,w)$ for any $(v,w)\in F$ we
have that no eigenvalue of $B^t B$ is greater than 1.
Therefore the smallest eigenvalue is at least $\gamma$.
This means that for all $(v,w)\in F$, $h_1 (v,v) \geq
\gamma \ h_2 (w,w)$.

Now if $(v^{\perp} ,w^{\perp} ) \in F^{\perp}$ we can find the
unique vector $ (v^{\perp} ,w^* )\in F$ whose first coordinate
is $v^{\perp}$. Then

$$h_1 (v^{\perp} ,v^{\perp}) =-h_2 (w^{\perp} ,w^* )
\leq \| w^{\perp} \| \|w^* \| \leq \frac{1}{\sqrt{\gamma}} \  \| w^{\perp} \|
\| v^{\perp} \| $$

\noindent
and so $h_2 (w^{\perp} ,w^{\perp} )\geq \gamma \ h_1 (v^{\perp}
,v^{\perp} )$

Therefore, if $h_1 (v,v)=1$ and
$(v,0)=f+f^{\perp}$ with $f=(I(v), w)\in F$ and
$f^{\perp} =(v^* ,-w)\in F^{\perp}$, we have that

$$h_1 (I(v),I(v) )\geq \gamma h_2 (w,w) \geq {\gamma}^2
h_1 (v^* ,v^* ).$$

\noindent
But $I(v)+v^* =v$, and therefore either $h_1 (v^* ,v^* )\geq 1/4$
or $h_1 (I(v),I(v) )\geq 1/4$. In any case,
$h_1 (I(v),I(v) )\geq (1/4) {\gamma}^2 $.

This means that ${\mu}_1 \geq (1/4){\gamma}^2 $ and so
$(\det I)^2 \geq 4^{-l} {\gamma}^{2l}$, proving (a).

\vspace{.5cm}

\noindent
{\bf b)} Given any $v\in V_1$ write $(v,0)=f_v +f_v^{\perp}$
(in $F\oplus F^{\perp}$). The map $L:V_1 \rightarrow
V_1 \oplus V_2$, $L(v)=f_v +\lambda f_v^{\perp}$, is
a monomorphism. Moreover, the image of $L$ is
included in the ${\bar{h}}_{\lambda}$-orthogonal complement
of $V_2$, $V_2^{{\perp}_{\lambda}}$,
and therefore $L$ gives an isomorphism between $V_1$
and $V_2^{{\perp}_{\lambda}}$.

Pick any $v\in V_1$. The map ${\pi}_1 \circ L
:V_1 \rightarrow V_1$ is an isomorphism. Therefore there
exists a unique $a\in V_1$ such that

$$v={\pi}_1 (f_a +\lambda f_a^{\perp} )= {\pi}_1
    \left( (f_{a1} ,f_{a2} )+\lambda (f_{a1}^{\perp} ,f_{a2}^{\perp} )
    \right) =f_{a1} +\lambda f_{a1}^{\perp} .$$

Since $(h_1 +h_2 )(f_a ,f_a^{\perp} )=0$, we have that

$$h_1 (f_{a1} ,f_{a1}^{\perp} ) +h_2 (f_{a2} ,f_{a2}^{\perp} )=0.$$

But since $f_a +f_a^{\perp} =(a,0)$, we have that $f_{a2} =-
f_{a2}^{\perp}$. Therefore

$$h_1 (f_{a1} ,f_{a1}^{\perp} )=h_2 (f_{a2} ,f_{a2} ) \geq 0.$$

There is a unique  $b\in V_1$ such that $(b,f_{a2}^{\perp} )
\in F$. Then since $h_1 (b,b) \leq h_2 (f_{a2}^{\perp},
f_{a2}^{\perp} )$ and $h_1 (b,f_{a1}^{\perp} ) +
h_2 (f_{a2}^{\perp} ,f_{a2}^{\perp} )=0$, we have that
$h_2 (f_{a2}^{\perp} ,f_{a2}^{\perp} ) \leq
h_1 (f_{a1}^{\perp} ,f_{a1}^{\perp} )$.

\vspace{.1cm}

Therefore

$$h_{\lambda} (v,v)={\bar{h}}_{\lambda} (f_a +\lambda f_a^{\perp} ,
f_a +\lambda f_a^{\perp} )=
\lambda (h_1 +h_2 ) (f_a ,f_a ) +{\lambda}^2 (h_1 +h_2 )
(f_a^{\perp} ,f_a^{\perp} )=$$

$$=\lambda h_1 (f_{a1} ,f_{a1} ) + \lambda h_2 (f_{a2} ,f_{a2} )
+{\lambda}^2 h_1 (f_{a1}^{\perp} ,f_{a1}^{\perp} ) + {\lambda}^2
h_2 (f_{a2}^{\perp} ,f_{a2}^{\perp} )$$

$$\leq h_1 (f_{a1} ,f_{a1} ) +\lambda h_1 (f_{a1} , f_{a1}^{\perp} )
+{\lambda}^2 h_1 (f_{a1}^{\perp} ,f_{a1}^{\perp} ) + {\lambda}
h_1 (f_{a1} ,f_{a1}^{\perp} ) = $$

$$=h_1 (f_{a1} +\lambda f_{a1}^{\perp} ,
f_{a1} +\lambda f_{a1}^{\perp} ) = h_1 (v,v).$$

Hence for any $v\in V_1$, $h_{\lambda} (v,v) \leq h_1 (v,v)$ and
(b) follows.

\end{proof}

We are now ready to prove our theorem.

\medskip
\noindent {\bf Theorem A.} {\it If the closed  manifold $M$ admits
an $\mathcal{F}$-structure then the minimal entropy of $M$ is 0. }

\medskip

\begin{proof}
Let $U_1 ,...,U_N$ be the open covering  corresponding to an
$\mathcal{F}$-structure on $M$;
with corresponding actions $._1 ,...,._N$ by tori
$T^{k_1},...,T^{k_N}$ on the coverings $\widetilde{U_i}$'s.

We can construct a {\it regular atlas} for the structure as
in \cite{Cheeger}, Lemma 1.2. Namely we construct a new open
cover $W_1 ,...,W_J$ of $M$ obtained by considering all
non-empty intersections of the $U_i$'s and then removing from
each set the ``unnecessary'' parts.
Each $W_i$ has a finite cover $\widetilde{W_i}$ where there
is defined an effective torus action.
For instance, if one had
$U_1 \cap U_2 \neq \emptyset$ then one would consider
$W_1 =U_1 \cap U_2$, $W_2 \subset U_1$, $W_3 \subset U_2$
so that $W_2 \cap W_3 =\emptyset$ and both are invariant through
the corresponding action (note that on $\widetilde{W_1}$ one has defined
a $T^{k_1 +k_2}$ action).

A Riemannian metric
$g$ on $M$ is called invariant
if on each of the open subsets $W_i$
of the $\mathcal{F}$-structure the corresponding sheaf of
torus acts by
isometries.
An invariant metric always exists,
at least after replacing the open subsets $W_i$
by slightly smaller ones.
Such a metric is constructed in
\cite{Cheeger}, Lemma 1.3.

Let us then fix a Riemannian metric $g$ on $M$ invariant
through the $\mathcal{F}$-structure. Each $W_i$ is,
essentially, the intersection of certain number of $U_i$'s.
Assume that $ W_1 ,...,W_J $ are ordered in a
non-increasing way with respect to the number of the
$U_i$'s intersecting.
Therefore, if $i>j$ and $W_i \cap W_j \neq \emptyset$
the torus action on $W_i$,
restricted to $W_i \cap W_j$, is embedded
in the action on $W_j$.
Consider smooth functions
$f_i :M \rightarrow [0,1]$, supported in $W_i$
which are constant along the orbits
and such that ${\{ f_i =1 \} }_{i=1,...,J}$ covers $M$.

Let $K=\sum_{i=1}^N k_i$ and let $\bar{M}=M\times T^K$ and
$\bar{g} =g+dx^2$ (where $dx^2$ is the standard Euclidean metric
on the $K$-torus). For each open subset $U_i$ consider the
following (free) $T^{k_i}$-action on $\widetilde{U_i} \times T^K$:

$$ ._{\bar{i}} : T^{k_i} \times
       \widetilde{U_i} \times T^K \rightarrow
       \widetilde{U_i} \times T^K$$

$$(\lambda ,(x,t_1 ,...,t_i ,...,t_N )) \mapsto
\left( \lambda ._i x,(t_1 ,...,\lambda t_i ,...,
t_N \right) $$

\noindent
where $t_j \in T^{k_j}$.

These formulas clearly define an $\mathcal{F}$-structure on $\bar{M}$. But
what is more important to us is that it is actually a polarized
$\mathcal{F}$-structure of positive rank.
Note that on the $W_i$'s all the torus actions corresponding to
the $U_i$'s which are intersecting glue together to get
a free torus action on $\widetilde{W_i } \times T^K$
(where the dimension of the torus acting is
the sum of the corresponding $k_i$'s).

Pull back the functions $f_i$ to obtain smooth functions $\bar{f_i}$
on $\bar{M}$. Note that the functions $\bar{f_i}$ are invariant
through both the torus action (on
$W_i \times T^K$) coming from $M$ and the canonical
$T^K$-action on the $T^K$-factor of $\bar{M}$. The same is true
for the metric $\bar{g}$.

Now we proceed to collapse the metric $\bar{g}$ along the orbits
of the $\mathcal{F}$-structure on $\bar{M}$. This is done in
\cite{Cheeger}, Theorem 3.1. We will describe the procedure, since
we need to make some computations on it. Fix a small $\delta >0$.
For technical reasons it is convenient to first replace $\bar{g}$
by ${\bar{g}}_0 =({\log}^2 \delta )\ \bar{g}$. We construct a
metric ${\bar{g}}_1$ on $\bar{M}$ by multiplying the metric
${\bar{g}}_0$ by ${\delta}^{\bar{f_1}}$ in the directions tangent
to the orbits of the torus action on $W_1 \times T^K$ (and leaving
the same metric in the directions orthogonal to the orbits). Note
that the $T^K$-action on $\bar{M}$ given by the canonical action
on the $T^K$-factor is isometric with respect to ${\bar{g}}_1$.
Repeating this procedure $J$-times we get a metric
${\bar{g}}_{\delta} ={\bar{g}}_J$ which is invariant under the
$T^K$-action just mentioned.

Let $g_{\delta}$ be the metric induced on $M=\bar{M} /T^K$. The
projection $(\bar{M} ,{\bar{g}}_{\delta} ) \rightarrow
(M,g_{\delta} )$ is a Riemannian submersion. Therefore the entropy
of $g_{\delta}$ is bounded above by the entropy of
${\bar{g}}_{\delta}$ (see Section 2.5). The entropy of
${\bar{g}}_{\delta}$ on the other hand is bounded above by
$(n-1)\sqrt{K_0 }$, where $K_0$ is an upper bound for the absolute
value of the sectional curvature of  ${\bar{g}}_{\delta}$ (see
Section 2.3). But it is proved in \cite{Cheeger}, Theorem 3.1,
that the sectional curvature of  ${\bar{g}}_{\delta}$ is bounded
independently of $\delta$.

Therefore we got that:

$${\rm h}_{{\rm top}}(g_{\delta}) \leq c_1 $$

\noindent
where $c_1$ is some constant independent of $\delta$.

We will now estimate the volume of $(M,g_{\delta})$.
We will do this by comparing the volume element of
$g_{\delta}$ with that of $g$.

Let $({\phi}_1 ,...,{\phi}_n )$ be a
$g$-orthonormal basis of $T_x M$.
Then

$${\rm dvol}(g_{\delta})= \sqrt{\det(g_{\delta}({\phi}_i ,{\phi}_j
))}{\rm dvol}(g).$$

Since the volume element at a point depends only on the value of
the metric at the point, it is the same to work on $W_i$ or on the
corresponding finite covering. Therefore from now on we will think
that we are working with a $\T$-structure to simplify the
notation. Fix any point $x\in M$ and any point $(x,t)\in \bar{M}$
which projects to $x$. We have to check how the volume element
changes at each step in the construction of ${\bar{g}}_{\delta}$.
Of course there is no change in the step $i$ if $x$ does not
belong to $W_i$. So let us assume for instance that $x\in W_1$.
Moreover, assume that $x$ is not a fixed point for the torus
action (the set of fixed points has volume 0 with respect to any
Riemannian metric). We want to compare the volume elements at $x$
of $g_1$ and $g$ ($g_i$ is of course the quotient of ${\bar{g}}_i$
under the $T^K$-action on $\bar{M}$).

Assume that the orbit through $x$ of the torus action
has dimension $l$. There is then an orthonormal set of vectors
${\omega}_1 ,...,{\omega}_l \in T_t (T^K )$ and some
linearly independent vectors
$v_1 ,...,v_l \in T_x M$
so that the vectors $(v_1 ,{\omega}_1 ),...,(v_l ,{\omega}_l )$
are tangent to the orbit on $\bar{M}$,
and the directions orthogonal to the ${\omega}_i$'s act trivially
on $M$ at $x$. Let $H$ be the subspace
of $T_{(x,t)} (M\times T^K )$ spanned by this $l$ vectors
(the tangent space to the orbit on $\bar{M}$),
let $V=<v_1 ,...,v_l > \ \subset T_x M$ be the tangent space to
the orbit in $M$ and let
$W=<w_1 ,...,w_l > \ \subset T_t (T^K )$.

Let $v_{l+1},...,v_n$ be a  $g$-orthonormal basis of the
space $g$-orthogonal to the orbit
(in $M$). Note that $v_{l+1},...,v_n$ are also
$g_1$-orthogonal to the orbit and $g_1 (v_{l+j} ,
v_{l+k} )= {\delta}_j^k (\log \delta )^2$. Therefore

$$\det \left( {\left( g_1 (v_i ,v_j \right) }_{1\leq i,j\leq n}
       \right) =  (\log \delta )^{2(n-l)}
  \det \left( {\left( g_1 (v_i ,v_j \right) }_{1\leq i,j\leq l}
        \right) .
$$

\vspace{.2cm}

Recall that the metric ${\bar{g}}_1$ is obtained by multiplying
by ${\delta}^{{\bar{f}}_1}$ the values of ${\bar{g}}_0$ on $H$.

From now on we restrict our attention to $V\oplus W$, since its
orthogonal complement plays no real role in the construction
of $g_1$.

We can assume that for any unitary tangent vector to any of the
tori (acting on any of the $W_i$),
the derivative of the action in that direction has $g$-norm
at most one.

Therefore we are under the hypothesis of our Linear Algebra
Lemma 6.1.

Consider now a $g$-orthonormal basis of $V$; call them $v^0_1 ,...,v^0_l $.
For each $i$ write

$$(v^0_i ,0)=h_i + h_i^{\perp} ,$$

\noindent
where $h_i \in H$ and $h_i ^{\perp}
 \in H^{\perp} $ (the $\bar{g}$-orthogonal
complement of $H$ in $V\oplus W$). Let $I(v^0_i )={\pi}_1 (h_i ) \in V$.

Now, for each $i=1,...,l$, consider the vector

$$w_i^0 = h_i + {\delta}^{{\bar{f}}_1} h_i^{\perp}.$$

\noindent
The vector $w_i^0$
is ${\bar{g}}_1$-orthogonal to the tangent space to the torus factor.
Its first coordinate is, of course, $I(v^0_i )+ {\delta}^{{\bar{f}}_1} \
{\pi}_1 (h_i^{\perp} )$.

Assume that  ${\bar{f}}_1 (x)=1$. Then

$$\det \left( g({\pi}_1 (w_i^0 ),{\pi}_1 (w_j^0 )) \right)
  =
  \det \left( g(I(v^0_i ),I(v^0_j )) \right)  \
+  \ o(\delta ).$$

Note also that:

$$g_1 ({\pi}_1 (w_i^0 ),
{\pi}_1 (w_j^0 )) =
{\bar{g}}_1 (w_i^0 ,w_j^0 )= {\log}^2 \delta \ \left(
{\delta}^{{\bar{f}}_1} \bar{g} (h_i
,h_j ) + {\delta}^{2{\bar{f}}_1} \ \bar{g} (h_i^{\perp} ,h_j^{\perp} )
\right) .$$

Therefore,

$$\frac{ \det \left( g_1 ({\pi}_1 (w_i^0 ),{\pi}_1 (w_j^0 )) \right) }{
        \det \left( g({\pi}_1 (w_i^0 ),{\pi}_1 (w_j^0 )) \right) }
= \frac{o({\log}^{2l} (\delta ) \ {\delta}^l )}{
det \left( g(I(v^0_i ),I(v^0_j )) \right) \ \ + \ \ o(\delta )}.$$

\vspace{.3cm}

Now, in the region where $\det \left( g(v_i ,v_j ) \right) >
{\delta}^{1/(4l)}$, we have from part {\bf (a)} of Lemma 6.1 that

$$\det \left( g(I(v^0_i ),I(v^0_j )) \right)
=(\det I)^2 \geq \frac{1}{4^l}
{\delta}^{1/2}$$

\noindent
and, therefore,

$${\rm dvol} (g_1 )= (\log \delta )^{n-l} \sqrt{ \frac{ \det
\left( g_1 ({\pi}_1 (w_i^0 ),{\pi}_1 (w_j^0 )) \right) }{
        \det \left( g({\pi}_1 (w_i^0 ),{\pi}_1 (w_j^0 )) \right) } }
     {\rm dvol}(g)=
o({\delta}^{1/4} {\log}^n \delta ){\rm dvol}(g).$$

Therefore the $g_1$-volume of the region
where $\det \left( g(v_i ,v_j ) \right) >
{\delta}^{1/(4l)}$ and $f_1 (x)=1$
approaches 0 as $\delta$ does.

The $g$-volume of the region $\det \left( g(v_i ,v_j ) \right)
<{\rho}^2$ is of the order of $\rho$. Therefore the
$g_1$-volume of the region where
$\det \left( g(v_i ,v_j ) \right) <{\delta}^{1/(4l)}$
is of the order of ${\log}^n (\delta ) {\delta}^{1/(8l)}$
(using part {\bf (b)} of Lemma 6.1)
and therefore it also approaches 0 with $\delta$.

This of course implies that ${\rm Vol}( \{ f_1 =1 \},g_{\delta})$
approaches 0 with $\delta$.

Finally note that in the passage from $\bar{g}$ to ${\bar{g}}_1$
there are two steps: first we multiply by ${\log}^2 \delta$ to
obtain $\bar{g_0}$ and then we collapse along the orbits
multiplying by ${\delta}^{\bar{f_1}}$. Lemma 6.1, part {\bf (b)},
tells us that the second of these steps does not
increase volumes (on $M$ with the quotient metric).
To go from ${\bar{g}}_1$ to ${\bar{g}}_2$ only the second step
is performed.
Therefore when passing from $g_1$ to
$g_2$ the volume of the region $f_1 =1$ will remain small,
while by taking $\delta$ small we can make the volume of the
region $f_2 =1$ small. Hence for $\delta$ small enough
the $g_{\delta}$-volume of the whole $M$ will be as
small as desired.

Since the entropy of $g_{\delta}$ is bounded above independently
of $\delta$, the theorem is proved.

\end{proof}

\section{Collapsing $\mathcal{F}$-structures and curvature invariants: Proof of Theorem B}
\label{yamabe}

There are many natural invariants of a smooth manifold which measure
the possible size of the curvature of a Riemannian metric of some
fixed volume. In this section we will recall some of them and study
what can be said about them for manifolds which admit
$\mathcal{F}$-structures. In every case we restrict attention
to the metrics verifying certain bounds on its curvature and
search for the infimum of the volumes.

Given a fixed closed smooth manifold $M$ we consider the following
subsets of the family $\mathcal{M}$ of all Riemannian metrics on
$M$:

$${\mathcal{M}}_{|K|} =\{ g: |K|\leq 1 \} $$

$${\mathcal{M}}_{K} = \{ g: K\geq -1 \} $$

$${\mathcal{M}}_{|r|} = \{ g: |r|\leq n-1 \} $$

$${\mathcal{M}}_r = \{ g : r \geq -(n-1) \} $$

$${\mathcal{M}}_{|s|} = \{ g : |s| \leq n(n-1) \} $$

$${\mathcal{M}}_s  = \{ g : s \geq -n(n-1) \} $$

\noindent
where $K$, $r$ and $ s$ denote as usual the sectional, Ricci and
scalar curvature, respectively. Now define (see
\cite{Gromov, LeBrun1} )

$${\rm MinVol} (M)= \inf_{g\in {\mathcal{M}}_{|K|} } {\rm
Vol}(M,g)$$

$${\rm Vol}_K (M) =\inf_{g\in {\mathcal{M}}_K } {\rm Vol}(M,g)$$

$${\rm Vol}_{|r|} (M) =\inf_{g\in {\mathcal{M}}_{|r|} } {\rm
Vol}(M,g)$$

$${\rm Vol}_r (M)=\inf_{g\in {\mathcal{M}}_r } {\rm Vol}(M,g)$$

$${\rm Vol}_{|s|} (M) =\inf_{g\in {\mathcal{s}}_{|K|} } {\rm
Vol}(M,g)$$

$${\rm Vol}_s (M) = \inf_{g\in {\mathcal{M}}_{|K|} } {\rm
Vol}(M,g)$$

Cheeger and Gromov proved that if $M$ is a closed manifold which
admits a polarized $\mathcal{F}$-structure of positive rank then
${\rm MinVol}(M)=0$.

It is easy to check the same proof of the theorem in the previous
section proves that if the closed  manifold $M$ admits an
$\mathcal{F}$-structure, then ${\rm Vol}_K (M)=0$. Of course, this
implies that ${\rm Vol}_r (M) = {\rm Vol}_s (M) = 0$.

More can be said about the scalar curvature.
Let us first recall some facts about the Yamabe invariant
(or {\it sigma constant} in \cite{Schoen}). More details
and references
can be found for instance in \cite{LeBrun2, Petean, Schoen}.

Given a conformal class of metrics $\mathcal{C}$ on $M$ the
Yamabe constant of $\mathcal{C}$, denoted by
$Y(M, \mathcal{C} )$, is the infimum of the integral of
the scalar curvature over all metrics in $\mathcal{C}$ of
volume 1 (integrating with respect to the volume element
of the same metric). The infimum is actually realized: this
is a very deep result obtained in several steps by
Yamabe, Trudinger, Aubin and Schoen. Metrics realizing
the infimum have constant scalar curvature and are usually
called Yamabe metrics. The Yamabe invariant of $M$, $Y(M)$, is
then defined as

$$Y(M)=\sup_{\mathcal{C}} Y(M,\mathcal{C} ).$$

$M$ admits a metric of strictly positive scalar curvature if and
only if $Y(M)>0$. In this case, if the dimension of $M$ is at
least 3,  $M$ also admits scalar flat metrics and so ${\rm
Vol}_{|s|} (M) ={\rm Vol}_s (M) =0$.

Now assume that $Y(M)\leq 0$. Let $g\in {\mathcal{M}}_s (M)$. Then
there exists a Riemannian metric $\hat{g} =e^f g$ conformal to $g$
with constant scalar curvature and with the same volume as $g$.
The scalar curvature of $\hat{g}$ can be written in terms of $f$
and $g$. From there it is easy to see that $s_{\bar{g}} \geq -1$
(see for instance \cite{Kobayashi}). But since $Y(M,
{\mathcal{C}}_g )\leq 0$, $s_{\bar{g}} \leq 0$. Therefore $\bar{g}
\in {\mathcal{M}}_{|s|} (M)$ and ${\rm Vol}(M,\bar{g})={\rm
Vol}(M,g)$. Hence:

\begin{Proposition} For any closed smooth manifold $M$ of dimension
greater than 2, ${\rm Vol}_{|s|} (M) ={\rm Vol}_s (M)$.
\end{Proposition}

Summarizing, we have proved the following:

\begin{Theorem} If $M$ admits an $\mathcal{F}$-structure,
$\mbox{\rm dim}\,M>2$,  then

$${\rm Vol}_K (M) ={\rm Vol}_r (M) ={\rm Vol}_{|s|} (M)={\rm
Vol}_s (M) =0.$$

\noindent
The last equality is equivalent to $Y(M)\geq 0$.
\end{Theorem}

Clearly this theorem implies Theorem B in the introduction.

 As we mentioned before there are plenty of examples of closed manifolds
$M$ which admit $\mathcal{F}$-structures and verify ${\rm
MinVol}(M)>0$. Also C. LeBrun proved (see \cite{LeBrun2, LeBrun3})
that, for instance, an elliptic compact complex surface collapses
with bounded Ricci curvature (i.e ${\rm Vol}_{|r|} =0$) if and
only if it is minimal. Therefore we have that, for instance, ${\rm
Vol}_{|r|}(T^4 \# {\overline\cp}^2 )>0$, while $T^4 \#
{\overline\cp}^2$ does admit an $\mathcal{F}$-structure.

\section{Minimal entropy in dimensions 4 and 5: Proofs of Theorems C, D and E}
\label{me45}

We will study in this section the minimal entropy of simply
connected manifolds of dimensions four and five. The aim is to
give an idea of up to what point the previous results can
be used to compute minimal entropies.

Let us begin with dimension four. Homeomorphism types have been
classified by Freedman. But the main feature in dimension four is
the comparison between homeomorphism classes and diffeomorphism
classes. Freedman's results say that the homeomorphism type of a
smooth closed simply connected 4-manifold is determined by the
intersection form. Not every possible intersection form can be
realized by a smooth manifold and the number of diffeotypes
corresponding to each homeotype is essentially unknown. With
regards to the question of which intersection forms are realized
by a smooth manifold it all comes down to the well-known
11/8-conjecture. Namely, the basic examples of (homeomorphism
types of) simply connected smooth four-manifolds are $S^4$, $S^2
\times S^2$, $\cp^2$ and $K3$. By taking connected sums of them
(with different orientations) one can realize many intersection
forms. Namely, connected sums of $\cp^2$'s realize all positive
definite intersection forms (by the well-known result of
Donaldson) and varying the orientations of some of the factors one
gets all odd forms. Finally, the complicated part of the analysis
comes from the indefinite even intersection forms. Let $H$ be the
intersection form of $S^2 \times S^2$:

$$H= \left(\begin{array}{cr}0&1\\1&0\end{array}\right)$$

\noindent
and let

$$E_8 =\left( \begin{array}{cccccccc}2&1&0&0&0&0&0&0\\
                               1&2&1&0&0&0&0&0\\
                               0&1&2&1&0&0&0&0\\
                               0&0&1&2&1&0&0&0\\
                               0&0&0&1&2&1&0&1\\
                               0&0&0&0&1&2&1&0\\
                               0&0&0&0&0&1&2&0\\
                               0&0&0&0&1&0&0&2
               \end{array}
       \right) .$$

Every even indefinite bilinear form is equivalent to
$kE_8 +lH$ for some integers $k$ and $l\geq 0$. Rohlin's
theorem says that for a smooth closed spin
4-manifold the signature is divisible by 16.
For simply connected 4-manifolds the spin condition
means exactly that the intersection form is even. Therefore we
have that for the intersection forms of smooth simply connected
4-manifolds, $k$ is even. The intersection form of the $K3$
surface is $-2E_8 +3H$. By taking connected sums of $K_3$'s
and $S^2 \times S^2$ we see that any such bilinear form can
be realized as the intersection form of a smooth 4-manifold if
$l\geq (3/2)|k|$. The $11/8$-conjecture says precisely that
these are exactly all the bilinear forms which come from
smooth simply connected 4-manifolds. Therefore in the
previous sections we have shown that:

\begin{Theorem} Assuming the $11/8$-conjecture, every closed simply
connected smooth 4-manifold is homeomorphic to one whose
minimal entropy is 0.
\end{Theorem}

There are simply connected compact complex surfaces of general
type which are homeomorphic to connected sums of $\cp^2$'s (with
different orientations). Nevertheless, they do not collapse with
bounded scalar curvature (see \cite{LeBrun2}) and so they cannot
admit $\mathcal{F}$-structures from the results of the previous
section. The following question therefore seems very intriguing:

\medskip

{\bf Question}: Is the minimal entropy of a simply connected
compact complex surface of general type positive?

\vspace{.5cm}

Let us now consider 5-dimensional manifolds. As we explained in
Section 3, closed simply connected smooth 5-manifolds have been
classified by S. Smale \cite{smale} and D. Barden \cite{Barden}.
We will use this classification to prove the next theorem which
clearly implies Theorem C in the introduction.

\begin{Theorem} Every simply connected closed smooth
5-manifold $M$ admits a $\T$-structure and hence ${\rm h}(M)=0$
and ${\rm Vol}_K (M) =0$. Moreover, $M$ admits a polarized
$\T$-structure and hence ${\rm MinVol} (M)=0$ unless $M$ is
cobordant to 0 and non-spin with $1<i(M)<\infty$.
\end{Theorem}

\begin{proof} We will prove that $M$ admits a
$\T$-structure and then apply Theorems A and B.

By Theorem 5.4 it is enough to show that each of the building
blocks of the classification (see Section 3) admits a
$\T$-structure.

Consider a smoothly embedded 2-sphere ${\bf S}$ representing
$j$-times a generator of $H_2 (S^2 \times S^3 , \Z)$ ($1<j<\infty
$). $M_j$ is obtained by performing surgery on ${\bf S}$. In the
same way $X_j$ is obtained by performing surgery on a sphere
representing $2^j$-times the generator of $H_2 (X_{\infty},\Z)$
(note that even multiples of the generator have trivial normal
bundles). This is easy to check since these manifolds are
characterized by their homology groups and whether they are spin
or not.

Of course, $X_0 =S^5$ and  $M_{\infty}=S^2 \times S^3$ admit free
$S^1$-actions. $X_{\infty}$ also admits a free $S^1$-action since
the Hopf action on $S^3$ commutes with the structure group of the
bundle. If we consider $X_{\infty}$ as the quotient of $S^3 \times
S^3 \subset \C^2 \times \C^2$ by the $S^1$-action

$$(w,(z_{1},z_{2},z_{3},z_{4}))\mapsto
(wz_{1},wz_{2},wz_{3},z_{4}) ,$$

\noindent then the Hopf action is given by complex multiplication
in the last two coordinates. But to construct $\T$-structures on
all the $X_j$'s consider the $S^1$-action on $X_{\infty}$ given by
complex multiplication on the last coordinate. This action has
fixed points, of course. Call this action $A_2$, and $A_1$ the
free ``Hopf''-action. The second homology of $X_{\infty}$ is
generated by the image (under the projection) of $\{ z_3 =0,z_4 =1
\}$. Call this 2-sphere ${\bf S _0}$. Now, given any small
$\varepsilon$ consider the 2-sphere

$$ {\bf S_{\varepsilon}}= \left\{ \left( z_1 ,z_2 , \varepsilon
z_2  , (1-{\varepsilon}^2 {\| z_2
    \|}^2 )^{1/2} \right) \right\} /S^1
\subset X_{\infty}.$$

\noindent ${\bf S_{\varepsilon}}$ is homologous to ${\bf S_0}$ and
they intersect only at $N=(1,0,0,1)$. If we set the imaginary part
of $z_4$ to be 0, we get the non-trivial $S^2$-bundle over $S^2$,
$S^3 \times S^2 /S^1 \subset S^3 \times S^3 /S^1$, which is
diffeomorphic to $\cp^2 \# {\overline\cp}^2$. We can modify ${\bf
S_0} \cup {\bf S_{\varepsilon}}$ inside $\cp^2 \#
{\overline\cp}^2$ to obtain a smoothly embedded sphere ${\bf S}$
representing twice the generator of $H_2 (X_{\infty},\Z)$. The
orbits of the $A_2$-action passing through ${\bf S}$ form an
embedded $S^2 \times S^1 \subset X_{\infty}$. Its normal bundle is
${\bf D} \times S^1$, where ${\bf D}$ is the $D^2$-bundle over
$S^2$ with Euler characteristic 4. ${\bf D}$ can be represented as
the quotient of $S^3 \times D^2$ under the $S^1$-action

$${\lambda} (z_1 ,z_2 ,z_3) =
  (\lambda z_1 ,\lambda z_2 ,{\lambda}^4 z_3 ).$$

\noindent There is then a canonical $S^1$-action $A_3$ on ${\bf
D}$ given by complex multiplication in the last coordinate. Define
a $\T$-structure on $X_{\infty}$ by leaving the $A_2$ action on
$X_{\infty} - {\bf D} \times S^1$ and giving to ${\bf D}\times
S^1$ the $A_3$ action. The zero section of ${\bf D}$ is the
embedded $S^2 \times S^1$ and is exactly the fixed point set of
$A_3$. For any $j$, $1< j <\infty$, consider a 2-sphere embedded
in $S^2 \times S^1$ representing $2^{j-1}$-times the generator of
the second homology group. This sphere represents $2^j$-times the
generator of $H_2 (X_{\infty},\Z)$. Its normal bundle is ${\bf D}
\times \re$, which is isomorphic to the trivial bundle $S^2 \times
D^3$. ${\bf D}$ is usually presented as the union of two copies of
$D^2 \times D^2$ glued along $S^1 \times D^2$ by the map

$$(\lambda ,z) \mapsto (\lambda , {\lambda}^4 z).$$

\noindent Namely, ${\lambda}^4$ is considered as a map $\gamma
:S^1 \rightarrow SO(2)$ and then we identify $(\lambda ,z)$ with
$(\lambda , \gamma (\lambda ) (z))$. The identification of ${\bf
D}\times \re$ with $S^2 \times D^3$ is obtained by an homotopy of
the loop $(\gamma ,1)$ in $SO(3)$ with the constant loop 1. The
action $A_3$ can then be viewed in $S^2 \times D^3$ as:

$${\lambda}.(x,y)=(x, {\lambda}_{\times} (\varphi (x)(y))),$$

\noindent for a map $\varphi :S^2 \rightarrow SO(3)$. Here
${\lambda}_{\times}$ means complex multiplication in the first two
(real) coordinates. Since ${\pi}_2 (SO(3))=1$ the map $\varphi$ is
null-homotopic. Therefore we can define an $S^1$-action on $S^2
\times (D_3 -\{ 0\} )$ which is equal to $A_3$ in an exterior
annulus and to

$${\lambda}.(x,y)=(x,{\lambda}_{\times} y)$$

\noindent in an inner annulus.

$X_j$ is obtained by deleting $S^2 \times D^3$ of $X_{\infty}$ and
gluing $D^3 \times S^2$ along the boundaries. Giving any
$S^1$-action to the $D^3$-factor of the glued $D^3 \times S^2$
clearly defines a $\T$-structure on $X_j$. These $\T$-structures
are not polarized.

\vspace{.4cm}

The Wu-manifold $X_{-1}= SU(3)/SO(3)$ admits a locally-free
$S^1$-action: simply embed $S^1$ in $SU(3)$ by sending $\lambda
\in S^1$ to the diagonal matrix with $\lambda$, $\lambda$,
${\lambda}^{-2}$ as the diagonal coefficients and then follow by
matrix multiplication.

\vspace{.4cm}

Finally, for any $j$, $1<j<\infty$, $M_j$ is obtained by
performing surgery on a sphere ${\bf S}$ representing $j$-times
the generator of $H_2 (S^2 \times S^3 ,\Z)$, which can be
represented by a completely transversal sphere for the Hopf action
on the $S^3$-factor (in the sense of Section 5). One then obtains
by Theorem \ref{surgery} a polarized $\T$-structure on $M_j$.

\vspace{.3cm}

This finishes the first part of the theorem. The last statement
follows because the fact that $M$ is either non-cobordant to zero
or it is cobordant to 0 but either $i(M) =0,1$ or $\infty$, means
that in the factorization of $M$ as connected sum of building
blocks only appear $M_j$'s, $X_{-1}$, $X_{1}$ and $X_{\infty}$ and
we have put polarized $\T$-structures on these manifolds
($X_{1}=X_{-1}\# X_{-1}$).

\end{proof}

\medskip

\noindent{\it Proof of Theorem D.} We shall make use of the
following remarkable fact which is a consequence of results M.
Gromov, Y. Yomdin and the Morse theory of the loop space. A proof
can be found in \cite{P}.

\begin{Theorem}Let $M$ be a closed simply connected smooth
manifold. Suppose that the loop space homology of $M$
\[\sum_{i=0}^{n}\mbox{\rm dim}\,H_{i}(\Omega M,k_{p})\]
grows exponentially with $n$ for some field of coefficients
$k_{p}$. Then, any $C^{\infty}$ Riemannian metric has positive
topological entropy. \label{poshtop}
\end{Theorem}

Let $M$ be a closed manifold obtained by taking connected sums of
copies of $S^{4}$, $\cp^{2}$, $\overline{\cp}^{2}$, $S^{2}\times
S^{2}$ and the $K3$ surface. Since $S^{4}$, $\cp^{2}$,
$\overline{\cp}^{2}$ and $S^{2}\times S^{2}$ admit a circle action
and the $K3$ surface admits a $\T$-structure by Theorem
\ref{telliptic}, it follows from Theorem \ref{connectedsums} that
$M$ admits a $\T$-structure. By Theorem A, the minimal entropy of
$M$ vanishes.

Suppose now that $M$ is diffeomorphic to one of the five manifolds
listed in Theorem D. By the results in Section \ref{zeroentropy}
each of these manifolds admits a smooth metric $g$ with $\h=0$ and
hence the minimal entropy problem can be solved for $M$.

On the other hand, suppose that the minimal entropy problem can be
solved for $M$. Since ${\rm h}(M)=0$ it follows that $M$ admits a
smooth metric with zero topological entropy. Theorem \ref{poshtop}
and Lemma \ref{lema2} imply that $M$ must be diffeomorphic to one
of the five manifolds listed in Theorem D. \qedsymbol

\medskip

\noindent{\it Proof of Theorem E.} Let $M$ be a closed simply
connected 5-manifold. Theorems C and A imply that the minimal
entropy of $M$ is zero.

Suppose now that $M$ is diffeomorphic to one of the four manifolds
listed in Theorem E. By the results in Section \ref{zeroentropy}
each of these manifolds admit a smooth metric $g$ with $\h=0$ and
hence the minimal entropy problem can be solved for $M$.

On the other hand, suppose that the minimal entropy problem can be
solved for $M$. Since ${\rm h}(M)=0$ it follows that $M$ admits a
smooth metric with zero topological entropy. Theorem \ref{poshtop}
and Corollary \ref{class5m} imply that $M$ must be diffeomorphic
to one of the four manifolds listed in Theorem E. \qedsymbol

\end{document}